\documentclass[11pt,namelimits,sumlimits,a4paper]{amsart}

\usepackage{amssymb}
\usepackage{amsmath}
\usepackage[mathscr]{eucal}
\usepackage{slashed}
\usepackage{amsthm}
\usepackage[latin1]{inputenc}
\usepackage[T1]{fontenc}
\usepackage[UKenglish]{babel}
\usepackage{amsfonts}
\usepackage{amscd}
\usepackage{latexsym}
\usepackage{cite}
\usepackage{comment}
\usepackage{times,psfrag}
\usepackage{mathtools}
\usepackage{bm}
\usepackage{esint}
\usepackage{layout}
\usepackage{lmodern}
\usepackage{cite}
\usepackage{comment}
\usepackage{color}
\usepackage{hyperref}

\numberwithin{equation}{section}
\newtheorem{thm}[equation]{Theorem}
\newtheorem{lemma}[equation]{Lemma}
\newtheorem{prop}[equation]{Proposition}
\newtheorem{cor}[equation]{Corollary}

\theoremstyle{definition}
\newtheorem{definition}[equation]{Definition}

\newtheorem{assumption}[equation]{Assumption}

\theoremstyle{remark}
\newtheorem{remark}[equation]{Remark}
\newtheorem*{remark*}{Remark}

\newcommand{\ie}{\emph{i.e.} }

\newcommand{\cf}{\emph{cf.} }

\newcommand{\vspan}{\operatorname{span}}

\newcommand{\beq}{\begin{equation}}
\newcommand{\eeq}{\end{equation}}
\newcommand{\bea}{\begin{eqnarray}}
\newcommand{\eea}{\end{eqnarray}}

\newcommand{\C}{\mathbb{C}}

\newcommand{\R}{\mathbb{R}}

\newcommand{\Z}{\mathbb{Z}}

\newcommand{\T}{\mathbb{T}}
\newcommand{\HH}{\mathbb{H}}

\newcommand{\Sph}{\mathbb{S}}
\newcommand{\ra}{\rightarrow}

\newcommand{\diag}{\operatorname{diag}}

\newcommand{\ad}{\textrm{ad}}
\newcommand{\Lie}[1]{\mathfrak{#1}}

\newcommand{\Ps}{\text{PS}}
\newcommand{\ext}{\text{ext}}
\newcommand{\inter}{\text{int}}
\newcommand{\Aut}{\text{Aut}}

\newcommand{\RS}{\R^{2} \times \Sph ^{1}}
\newcommand{\Ann}{\text{Ann}}

\def\co{\colon\thinspace}

\begin{document}
\raggedbottom

\title{A gluing construction for periodic monopoles}

\author[L.~Foscolo]{Lorenzo~Foscolo}
\address{Mathematics Department, State University of New York at Stony Brook}
\email{lorenzo.foscolo@stonybrook.edu}

\date{\today}


\begin{abstract}
In \cite{Cherkis:Kapustin:1} and \cite{Cherkis:Kapustin:2} Cherkis and Kapustin introduced the study of periodic monopoles (with singularities), \ie monopoles on $\R ^{2} \times \Sph ^{1}$ possibly singular at a finite collection of points. Four-dimensional moduli spaces of periodic monopoles with singularities are expected to provide examples of gravitational instantons, \ie complete hyperk\"ahler $4$--manifolds with finite energy.

In a previous paper \cite{Foscolo:Deformation} we proved that the moduli space of charge $k$ periodic monopoles with $n$ singularities is either empty or generically a smooth hyperk\"ahler manifold of dimension $4(k-1)$. In this paper we settle the existence question, constructing periodic monopoles (with singularities) by gluing methods.
\end{abstract}

\maketitle

\section{Introduction}\label{sec:Introduction}

Let $(X,g)$ be an oriented Riemannian $3$--manifold and $P \ra X$ a principal $G$--bundle, where $G$ is a compact Lie group. (Magnetic) \emph{monopoles} are solutions $(A,\Phi )$ to the \emph{Bogomolny equation}
\begin{equation}\label{eqn:Bogomolny}
\ast F_{A}=d_{A}\Phi.
\end{equation}
Here $\ast$ is the Hodge star operator of $(X,g)$; $F_{A}$ is the curvature of a connection $A$ on $P$; and $\Phi$, the \emph{Higgs field}, is a section of the adjoint bundle $\ad (P)$. The moduli space of monopoles on $P \ra X$ is the space of equivalence classes of solutions to \eqref{eqn:Bogomolny} with respect to the action of the gauge group $\Aut (P)$. The Bogomolny equation is the dimensional reduction of the Yang-Mills anti-self-duality (ASD) equation, \ie monopoles on $X$ are circle-invariant instantons on $X \times \Sph ^{1}$.

An immediate consequence of equation (\ref{eqn:Bogomolny}) and the Bianchi identity is
\begin{equation}\label{eqn:Harmonic:Higgs:Field}
d_{A}^{\ast}d_{A}\Phi =0.
\end{equation}
In particular, when $X$ is compact smooth monopoles coincide with reducible (assuming $\Phi \neq 0$) flat connections. In order to find irreducible solutions to (\ref{eqn:Bogomolny}) one has to consider a non-compact base manifold $X$, in the sense that either $X$ is complete or we allow for singularities of the fields $(A,\Phi)$, or a combination of the two possibilities, as in this paper.

The classical case of smooth monopoles on $\R ^{3}$ has been investigated from many different points of view, \cf Atiyah and Hitchin's book \cite{Atiyah:Hitchin}. An important property of the moduli spaces of monopoles on $\R ^{3}$ is that they are hyperk\"aler manifolds by virtue of an infinite dimensional hyperk\"ahler quotient. In the 1980's Atiyah and Hitchin found an explicit formula for the metric on the moduli space of centred charge $2$ $SU(2)$ monopoles on $\R ^{3}$. From this formula it follows that the Atiyah--Hitchin manifold is a gravitational instanton (a complete hyperk\"ahler $4$--manifold with finite $L^{2}$--norm of the curvature) of type ALF: the volume of large geodesic balls of radius $r$ grows like $r^{3}$, the complement of a large ball is a circle bundle over $\R ^{3}/\Z_{2}$ and the metric is asymptotically adapted to this circle fibration.

Pursuing the idea that moduli spaces of solutions to dimensional reductions of the Yang--Mills ASD equations on $\R ^{4}$ are ``a natural place to look for gravitational instantons'' \cite{Cherkis:Talk}, in \cite{Cherkis:Kapustin:1}, \cite{Cherkis:Kapustin:3} and \cite{Cherkis:Kapustin:2} Cherkis and Kapustin introduced the study of \emph{periodic monopoles}, \ie monopoles on $\R ^{2} \times \Sph ^{1}$, possibly with isolated singularities at a finite collection of points. They argued that, when $4$--dimensional, moduli spaces of periodic monopoles (with singularities) are gravitational instantons of type ALG: the volume of large balls grows quadratically and the metric is asymptotically adapted to a fibration by $2$--dimensional tori.

In \cite{Foscolo:Deformation} we proved that for generic choices of parameters the moduli space $\mathcal{M}_{n,k}$ of $SO(3)$ periodic monopoles of charge $k$ with $n$ singularities is either empty or a smooth hyperk\"ahler manifold of dimension $4(k-1)$. Here the choice of generic parameters guarantees that $\mathcal{M}_{n,k}$ does not contain reducible solutions. In this paper we address the existence question and construct periodic solutions to \eqref{eqn:Bogomolny} by gluing methods. The main result of the paper is the following theorem. We refer to Corollary \ref{cor:Existence} for a more precise statement.

\begin{thm}\label{thm:Main:Theorem}
Under appropriate assumptions on the parameters defining the boundary conditions, there exist irreducible $SO(3)$ periodic monopoles (with singularities), \ie the moduli space $\mathcal{M}_{n,k}$ contains smooth points.
\end{thm}

Monopoles on $\R ^{3}$ (with structure group $SU(2)$ and without singularities) were themselves constructed via gluing methods in a seminal work by Taubes \cite[Theorem 1.1 \S IV.1]{Jaffe:Taubes}. Furthermore, Cherkis and Kapustin's physically-motivated computation of the asymptotics of the $L^{2}$--metric on the moduli space of periodic monopoles \cite{Cherkis:Kapustin:3} is based on the idea that a charge $k$ monopole is a non-linear superposition of particle-like charge $1$ components. The theorem confirms this expectation and we plan to exploit the gluing construction to recover Cherkis and Kapustin's asymptotic formula for the $L^{2}$--metric on the moduli spaces in a future paper.

The main steps and ingredients of Taubes's original gluing construction for Euclidean monopoles without singularities are:
\begin{itemize}
\item[(i)] Charge 1 monopoles on $\R ^{3}$ are completely explicit, as we will see in Section \ref{sec:Prasad:Sommerfield:Monopole}. Up to translations and scaling there exists a unique solution, localised around the origin in $\R^{3}$. 
\item[(ii)] Given $k$ points far apart in $\R^{3}$, Taubes constructs an approximate solution to (\ref{eqn:Bogomolny}) patching together $k$ charge 1 monopoles each localised around one of these points; the choice of gluing maps accounts for a further $k-1$ parameters.
\item[(iii)] The approximate solution is deformed to a genuine monopole by an application of the implicit function theorem.
\end{itemize}

The first difficulty to implement the construction in the periodic case is that not even charge $1$ periodic monopoles are explicitly known. In fact, numerical experiments of Ward \cite{Ward} show that a very different behaviour should be expected depending on the sign of the \emph{mass}, the constant term $v$ in the expansion of $|\Phi|$ at infinity, \cf Definition \ref{def:Boundary:Conditions}. When $v$ is positive and large, charge $1$ periodic monopoles are concentrated in an almost spherical region around their centre. When the mass is negative and large in absolute value, the monopoles are instead localised in a slab containing \emph{two} maxima of the energy density.

As a consequence, the construction of a charge $k$ periodic monopole as a superposition of $k$ charge $1$ monopoles can be carried out only when the charge $1$ constituents have large positive mass. There are two ways of arranging this. On one side one can consider periodic monopoles with large mass $v$. By scaling, the large mass limit $v \ra +\infty$ is equivalent to the large radius limit $\R ^{2} \times \R /2\pi v\Z \ra \R ^{3}$. Here nothing is special about the case $X=\R ^{2} \times \Sph ^{1}$ and it is conceivable that large mass monopoles exist on any $3$--manifold satisfying appropriate conditions. More interestingly, we will exploit the fact that the Green's function of $\R ^{2} \times \Sph ^{1}$ grows logarithmically at infinity (\cf Lemma \ref{lem:Asymptotics:Periodic:Dirac:Higgs:Field}) to construct periodic monopoles with arbitrary mass $v$ and $n<2(k-1)$ singularities. These solutions are described qualitatively as the superposition of widely separated charge $1$ components which get more and more concentrated around their respective centres as these recede from each other.

\subsection*{Plan of the paper}
The proof of Theorem \ref{thm:Main:Theorem} is articulated into various steps. The basic building blocks used in the gluing construction---periodic Dirac monopoles and charge $1$ Euclidean monopoles---are introduced in Section \ref{sec:Preliminaries} together with further background material.

The starting point of the construction is a certain singular solution to the Bogomolny equation described in Section \ref{sec:Sum:Dirac:Monopoles}: given singularities $p_{1}, \ldots, p_{n}$ and $k$ additional well-separated points $q_{1}, \ldots , q_{k}$, we construct a reducible solution to the Bogomolny equation on $(\RS) \setminus \{ p_{1}, \ldots, p_{n}, q_{1}, \ldots, q_{k} \}$ by taking a sum of periodic Dirac monopoles. Consideration of the asymptotic behaviour of this reducible solution leads to the definition of boundary conditions as in Cherkis--Kapustin \cite{Cherkis:Kapustin:1} and \cite{Cherkis:Kapustin:2}.

We want to think of $q_{1}, \ldots, q_{k}$ as the centres of highly concentrated charge $1$ monopoles. To this end it is necessary to assume that either
\begin{itemize}
\item[(A)] the mass $v$ is sufficiently large, or
\item[(B)] when the number of singularities $n$ is less than $2(k-1)$, $q_{1}, \ldots , q_{k}$ are sufficiently far away from each other and from the singularities $p_{1}, \ldots, p_{n}$.
\end{itemize}
We refer to (A) and (B) as the \emph{high mass} and \emph{large distance} case respectively. 
Under either of these hypothesis, in Section \ref{sec:Approximate:Solutions} we construct initial approximate solutions to the Bogomolny equation by gluing scaled Euclidean charge $1$ monopoles in a neighbourhood of $q_{1}, \ldots, q_{k}$ to resolve the singularities of the sum of periodic Dirac monopoles. By varying the centres and phases (thought of as fixing the choice of gluing maps) of the glued-in charge $1$ monopoles, we obtain a $4(k-1)$--dimensional family of inequivalent approximate solutions.

The next step of the construction is to deform the initial approximate solutions into genuine monopoles by means of the Implicit Function Theorem. The crucial step, tackled in Section \ref{sec:Linear}, is to study the linearised equation. A first difficulty arises from the fact that, if one fixes the boundary conditions (\ie works with weighted Sobolev spaces forcing certain decay), there is a $3$--dimensional space of obstructions to the solvability of the linearised equation. There are two ways to proceed:
\begin{itemize}
\item[(i)] enlarge the Banach spaces in which to solve the Bogomolny equation by allowing the appropriate changes of asymptotics;
\item[(ii)] consider the centre of mass of the centres of the glued-in Euclidean charge $1$ monopoles as a free parameter to be fixed only at the end of the construction to compensate for the obstructions.
\end{itemize}
We follow this second approach as it seems more appropriate to construct a whole $4(k-1)$--parameter family of monopoles in a fixed moduli space. In order to study the linearised problem, we adopt the strategy of studying the linearised equation separately on the building blocks, the charge $1$ Euclidean monopoles and the sum of periodic Dirac monopoles. In the former case, there are no obstructions to the solvability of the linearised equation and the use of weighted Sobolev spaces allows to obtain uniform estimates for the norm of a right inverse. In the latter case, we can solve the linearised equation in the chosen weighted Sobolev spaces only modulo obstructions. Furthermore, for technical reasons we have to distinguish between the high mass and large distance case (A) and (B) above.
\begin{itemize}
\item[(A)] Under the assumption that the points $q_{1}, \ldots, q_{k}$ are contained in a fixed compact set of $\RS$ and that the mass $v$ is sufficiently large, we use the weighted Sobolev spaces and estimates of \cite{Foscolo:Deformation} without major difficulties.
\item[(B)] When the points $q_{1}, \ldots, q_{k}$ move off to infinity, instead, an additional technical difficulty arises from the following fact: it is well-known that for all $f \in C^{\infty}_{0}(\R ^{2})$ with mean value zero there exists a bounded solution $u$, unique up to the addition of a constant, to $\triangle u =f$ with $\| \nabla u \| _{L^{2}} < \infty$. However, if $f$ is supported on the union of balls $B_{1}(z_{1}) \cup B_{1}(z_{2})$, say, with non-zero mean value on each of them, then $\| \nabla u \| ^{2}_{L^{2}} \geq C \log{|z_{1}-z_{2}|}$. As a consequence, in the large distance case the linearised operator is not well-behaved on a certain finite-dimensional space of sections which has to be dealt with separately.  
\end{itemize}
In Section \ref{sec:Linearised:Equation:Modulo:Obstructions} we patch together the local inverses of the linearised operator on the different building blocks and by a simple iteration solve the linearised equation globally modulo obstructions. Finally, in Section \ref{sec:Deformation} we conclude the construction of solutions to the Bogomolny equation satisfying the required boundary conditions by an application of the Implicit Function Theorem.

\subsection*{Aknowledgments}
The results of this paper are part of the author's Ph.D. thesis at Imperial College London. He wishes to thank his supervisor Mark Haskins for his support. Olivier Biquard guided early stages of this project; we thank him for suggesting us this problem. The author is grateful to Hans-Joachim Hein for discussions while this work was being developed and to Simon Donaldson and Michael Singer for their careful comments on an early version of this paper. The writing process was completed while the author was a Simons Instructor at Stony Brook University.

\section{Preliminaries}\label{sec:Preliminaries}

In this section we describe in some details the simple components to be patched together in the gluing construction. We begin with periodic Dirac monopoles, \ie abelian solutions to the Bogomolny equation on $\R ^{2} \times \Sph ^{1}$ with a singularity at one point, recalling the asymptotic expansions proved in \cite{Foscolo:Deformation}. Secondly, we collect the main properties of the basic Euclidean monopole, the Prasad--Sommerfield monopole. As a preliminary and mainly to fix some notation we give a brief overview of the deformation theory of monopoles on an arbitrary $3$--manifold.

\subsection{Deformation theory of monopoles}

Let $(X,g)$ be a non-compact oriented Riemannian $3$--manifold endowed with a principal $G$--bundle $P \ra X$.
Denote by $\mathcal{C}$ the infinite dimensional space of smooth pairs $c=(A,\Phi)$, where $A$ is a connection on $P \ra X$ and $\Phi \in \Omega ^{0}(X;\ad P)$ a Higgs field. Since $X$ is non compact, elements of $\mathcal{C}$ have to satisfy appropriate boundary conditions, which we suppose to be included in the definition of $\mathcal{C}$. The configuration space $\mathcal{C}$ is an affine space; the underlying vector space is the space of sections
\[
\Omega (X;\ad P):=\Omega ^{1}(X;\ad P)\oplus\Omega ^{0}(X;\ad P)
\]
satisfying appropriate decay conditions.

Let $\mathcal{G}$ be the group of bounded smooth sections of $\Aut(P)$ which preserve the chosen boundary conditions. Here $g \in Aut(P)$ acts on a pair $c=(A,\Phi) \in \mathcal{C}$ by $c \mapsto c + (d_{1}g)g^{-1}$, where\begin{equation}\label{eqn:Linearisation:Gauge:Action}
d_{1}g = -\left( d_{A}g, [\Phi,g] \right) \in \Omega (X; \ad P).
\end{equation}

Consider the gauge-equivariant map $\Psi\co \mathcal{C} \ra \Omega ^{1}(X;\ad P)$ defined by
\[
\Psi (A,\Phi )= \ast F_{A}-d_{A}\Phi.
\]
By fixing a base point $c=(A,\Phi )\in\mathcal{C}$ we write $\Psi (A+a,\Phi +\psi )=\Psi (c)+d_{2}(a,\psi )+(a,\psi )\cdot (a, \psi )$ for all $(a,\psi )\in\Omega (X; \ad P)$. The linearisation $d_{2}$ of $\Psi$ at $c$ and the quadratic term are defined by:
\begin{equation}\label{eqn:Linearisation:Bogomolny}
d_{2}(a,\psi )=\ast d_{A}a-d_{A}\psi +[\Phi ,a]
\end{equation}
\begin{equation}\label{eqn:Quadratic:Term:Bogomolny}
(a,\psi )\cdot (a, \psi )=\ast [a,a]-[a,\psi ]
\end{equation}

The linearisation at $c$ of the action of $\mathcal{G}$ on $\mathcal{C}$ is the operator $d_{1}\co \Omega ^{0}(X;\ad\,P) \ra \Omega (X; \ad P)$ defined as in \eqref{eqn:Linearisation:Gauge:Action}. We couple $d_{2}$ with $d_{1}^{\ast}$ to obtain an elliptic operator
\begin{equation}\label{eqn:Dirac:Operator}
D=d_{2}\oplus d_{1}^{\ast }\co \Omega (X; \ad P)\longrightarrow \Omega (X; \ad P).
\end{equation}
The moduli space $\mathcal{M}$ of monopoles in $\mathcal{C}$ is $\mathcal{M}=\Psi ^{-1}(0)/\mathcal{G}$. If $\mathcal{M}$ is smooth, the tangent space $T_{[c]}\mathcal{M}$ at a point $c=(A,\Phi)$ is identified with $\ker D$.

The operator $D$ is a twisted Dirac operator on $\Omega (X; \ad P)$. To see this, recall that Clifford multiplication is defined by
\begin{equation}\label{eqn:Clifford:Multiplication}
\gamma (\alpha)\beta = \alpha \wedge \beta - \alpha ^{\sharp} \lrcorner\, \beta
\end{equation}
for a $1$--form $\alpha$ and a $k$--form $\beta$. Then $D = \tau\slashed{D}_{A} + [\Phi ,\cdot\, ]$, where $\slashed{D}_{A}$ is the twisted Dirac operator
\[
\Omega ^{1} \oplus \Omega ^{0} \xrightarrow{(\text{id},\ast)} \Omega ^{1} \oplus \Omega ^{3} \xrightarrow{\gamma\, \circ \nabla _{A}} \Omega ^{2} \oplus \Omega ^{0} \xrightarrow{(\ast, \text{id})} \Omega ^{1} \oplus \Omega ^{0}.
\]
and $\tau$ is a sign operator with $\tau = 1$ on $1$--forms and $\tau =-1$ on $0$--forms.

\subsection{Periodic Dirac monopole}\label{sec:Periodic:Dirac:Monopole}

When the structure group is $G=U(1)$, the Bogomolny equation \eqref{eqn:Bogomolny} reduces to a linear equation. By \eqref{eqn:Harmonic:Higgs:Field} the Higgs field $\Phi$ is a harmonic function such that $\frac{\ast d\Phi}{2\pi i}$ represents the first Chern class of a line bundle. The assumption that $\Phi$ is bounded thus forces $\Phi$ to be constant. Non-trivial abelian solutions, so called \emph{Dirac monopoles}, are obtained by allowing an isolated singularity. On $\R^{3}$ these are defined as follows.

\begin{definition}\label{def:Euclidean:Dirac:Monopole}
Fix a point $q \in \R ^{3}$ and let $H_{q}$ denote the radial extension of the Hopf line bundle to $\R ^{3} \setminus \{ q \}$. Fix $k \in \Z$ and $v \in \R$. The \emph{Euclidean Dirac monopole} of charge $k$ and mass $v$ with singularity at $q$ is the abelian monopole $(A^{0}, \Phi^{0})$ on $H_{q}^{k}$, where
\[
\Phi^{0} =i\left( v-\frac{k}{2|x-q|} \right),
\]
$x \in \R ^{3}$, and $A^{0}$ is the $SO(3)$--invariant connection on $H_{q}^{k}$ with curvature $\ast d\Phi$.
\end{definition}

Periodic Dirac monopoles are defined in a similar way. Fix coordinates $(z,t)\in \C \times \R /2\pi \Z$ and a point $q=(z_{0},t_{0})\in \R ^{2} \times \Sph ^{1}$. Line bundles of a fixed degree on $(\RS) \setminus \{ q \}$ differ by tensoring by flat line bundles. We can distinguish connections with the same curvature by comparing their holonomy around loops $\gamma _{z}:=\{ z \} \times \Sph ^{1}_{t}$ for $z \neq z_{0}$. Set $\theta _{q}=\text{arg}(z-z_{0})$ and fix an origin in the circle parametrised by $\theta _{q}$. Denote by $L_{q}$ the degree $1$ line bundle on $(\RS) \setminus \{ q \}$ with connection $A_{q}$ whose holonomy around $\gamma _{z}$ is $e^{-i\theta _{q}}$. Any line bundle of degree $1$ is of the form $L_{q} \otimes L_{b}$ for some flat line bundle $L_{b}$.

\begin{definition}\label{def:Periodic:Dirac:Monopole}
Fix a point $q \in \RS$. The \emph{periodic Dirac monopole} of charge $k \in \Z$, with singularity at $q$ and twisted by the flat line bundle $L_{v,b}$ for some $v \in \R$ and $b \in \R /\Z$ is the pair $(A, \Phi)$ on $L_{q}^{k} \otimes L_{v,b}$, where
\[
\Phi =i \left( v+kG_{q} \right)
\]
and up to gauge transformations the connection is $A= kA_{q}+ib\, dt$. Here $G_{q}$ is the Green's function of $\RS$ with singularity at $q$ defined in Lemma \ref{lem:Asymptotics:Periodic:Dirac:Higgs:Field} below.
\end{definition}

In \cite{Foscolo:Deformation} we derived asymptotic expansions for the Green's function $G_{q}$ and the connection $A_{q}$, both at infinity and close to the singularity. As these expansions will be essential for the gluing construction, we recall them here. By taking coordinates centred at $q \in \RS$, we can assume that the singularity is located at $q=0$. We use polar coordinates $z=re^{i\theta} \in \C$.

\begin{lemma}[{Lemma 3.4 in \cite{Foscolo:Deformation}}]\label{lem:Asymptotics:Periodic:Dirac:Higgs:Field}
There exists a Green's function $G$ of $\RS$ with singularity at $0$ such that the following holds.
\begin{itemize}
\item[(i)] There exists a constant $C_{1}>0$ such that
\[
\left| \nabla ^{k}\left( G(z,t)-\frac{1}{2\pi }\log {r}\right)\right| \leq C_{1}e^{-r}
\]
for all $r\geq 2$ and $k=0,1,2$.
\item[(ii)] There exists a constant $C_{2}>0$ such that 
\[
\left| \nabla ^{k}\left( G(z,t)-\frac{a_{0}}{2}+\frac{1}{2\rho }\right)\right| \leq C_{2}\rho ^{2-k}
\]
for all $(z,t)$ with $\rho=\sqrt{r^{2}+t^{2}}< \frac{\pi}{2}$ and $k=0,1,2$.
\end{itemize}
\end{lemma}

Fix a constant $v \in \mathbb{R}$ and consider the Higgs field $\Phi =iv+iG$. The $2$--form $i\ast dG$ represents the curvature of a line bundle $L$ over $(\RS)\setminus\{ 0\}$. In a neighbourhood of the singularity $L$ is isomorphic to the Hopf line bundle extended radially from a small sphere $\Sph ^{2}$ enclosing the origin. Any connection $A$ on $L$ with $F_{A}=\ast d\Phi$ is asymptotically gauge equivalent to $A^{0}$ as $\rho \ra 0$.

At infinity $L$ is isomorphic to the radial extension of a line bundle of degree $1$ over the torus $\mathbb{T}^{2}_{\infty}$. We choose a representative for $A$ in this asymptotic model as follows. Consider the connection $A^{\infty}=-i\frac{t}{2\pi}d\theta$ on the trivial line bundle $\underline{\mathbb{C}}$ over $\Sph ^{1}_{\theta} \times \mathbb{R}_{t}$. Let $\tau$ be the map of $\underline{\C}$ into itself defined by $\tau (e^{i\theta },t,\xi )=(e^{i\theta },t+2\pi,e^{i\theta}\xi )$ and define a line bundle with connection over $\mathbb{T}^{2}_{\theta ,t}$ as the quotient $(\underline{\mathbb{C}},A^{\infty })/\tau$. As $r\rightarrow\infty$, up to gauge transformations $A$ is asymptotic to $A^{\infty}+i\alpha\, d\theta +ib\, dt$ for some $\alpha ,b\in\mathbb{R}/\mathbb{Z}$. The monodromy of this limiting connection is $e^{-i\theta-2\pi ib}$ around the circle $\{\theta\}\times \Sph ^{1}_{t}$ and $e^{it-2\pi i\alpha}$ around the circle $\Sph ^{1}_{\theta}\times\{ t\}$. While $b$ can be chosen arbitrarily, by \cite[Lemma 3.6]{Foscolo:Deformation} $\alpha$ is fixed by the Bogomolny equation \eqref{eqn:Bogomolny}, $\alpha =\frac{1}{2}$ modulo $\Z$.

\begin{lemma}[{Lemma 3.6 in \cite{Foscolo:Deformation}}]\label{lem:Asymptotics:Periodic:Dirac:Connection}
Fix parameters $(v,b) \in \R \times \R /\Z$. Let $(A,\Phi )$ be a solution to (\ref{eqn:Bogomolny}) with $\Phi =i \left( v+G \right)$ and such that the holonomy of $A$ around circles $\{ re^{i\theta} \} \times \Sph ^{1}_{t}$, $r \neq 0$, is $e^{-i\theta-2\pi ib}$.
\begin{itemize}
\item[(i)] In the region where $r\geq 2$ the connection $A$ is gauge equivalent to
\[
A^{\infty}+\frac{i}{2}\, d\theta +ib\, dt+a
\]
for a $1$--form $a$ such that $d^{\ast}a=0=\partial _{r}  \lrcorner\, a$ and $|a|+|\nabla a|=O(e^{-r})$.
\item[(ii)] In a ball of radius $\frac{\pi}{2}$ centred at the singular point $z=0=t$, $A$ is gauge equivalent to $A^{0}+a'$ where $|a'|+\rho |\nabla a'|=O(\rho ^{2})$ and $d^{\ast}a'=0=\partial _{\rho}\lrcorner \, a'$.
\end{itemize}
\end{lemma}

Given an arbitrary point $q=(z_{0},t_{0})$ in $\RS$ the same formulas describe the asymptotic behaviour of the periodic Dirac monopole $(A_{q},\Phi _{q})$ with singularity at $q$ in coordinates centred at $q$. It will be useful to express the behaviour of $(A_{q},\Phi _{q})$ at large distances from $q$ in a fixed coordinate system.

\begin{lemma}[{Lemma 3.7 in \cite{Foscolo:Deformation}}]\label{lem:Asymptotics:Periodic:Dirac:Translations}
For $r\geq 2|z_{0}|$ we have
\[
\begin{gathered}
\frac{1}{i}\Phi _{q}(z,t) = v+\frac{1}{2\pi}\log{r}-\frac{1}{2\pi}\, \text{Re}\left( \frac{z_{0}}{z}\right) +O(r^{-2}) \\
A_{q}(z,t) = A^{\infty } + ib\, dt + i\, \frac{t_{0}+\pi}{2\pi}\, d\theta -\frac{i}{2\pi}\, \text{Im}\left( \frac{z_{0}}{z}\right) dt+O(r^{-2}).
\end{gathered}
\]
\end{lemma}

Finally, notice that the parameters $(v,b) \in \R \times \R /\Z$ are related to rotations and dilations. By a rotation in the $z$--plane, we can always assume that $b=0$. On the other hand, given any $\lambda>0$ consider the homothety
\[
h_{\lambda }:\R ^{2} \times \R/2\pi\Z\longrightarrow \R ^{2} \times \R/2\pi\lambda\Z
\]
of ratio $\lambda$. Since the Bogomolny equation is the dimensional reduction of the ASD equation, which is conformally invariant, $(h_{\lambda }^{\ast}A ,\lambda\, h_{\lambda }^{\ast }\Phi )$ is a monopole on $\R ^{2} \times \R/2\pi\Z$ if and only if $(A,\Phi )$ solves the Bogomolny equation on $\R ^{2} \times \R/2\pi\lambda\Z$. A simple but crucial observation for the gluing construction is that given a periodic Dirac monopole $(A_{q},\Phi_{q})$ with mass $v$, then as $v\rightarrow \infty$
\[
\lambda^{-1}h_{\lambda ^{-1}}^{\ast }\Phi\longrightarrow i\left( 1-\frac{1}{2\sqrt{r^{2}+t^{2}}}\right),
\]
where we set $\lambda =v+\frac{a_{0}}{2}$. In other words, the limit $v\rightarrow \infty$ corresponds to the limit $\mathbb{R}^{2}\times \Sph ^{1}\rightarrow\mathbb{R}^{3}$ and in this limit a periodic Dirac monopole converges to an Euclidean Dirac monopole.

\subsection{Charge 1 monopoles on $\mathbb{R}^{3}$}\label{sec:Prasad:Sommerfield:Monopole}

As we will see in the next section, periodic Dirac monopoles will serve us to construct an initial \emph{singular} solution to the Bogomolny equation. We now describe the model solution that we will use to desingularise this initial solution. In 1975 Prasad and Sommerfield \cite{Prasad:Sommerfield} found an explicit smooth finite energy solution to the Bogomolny equation on $\R ^{3}$ with structure group $SU(2)$. By translations and scaling, this explicit solution accounts for all $SU(2)$ Euclidean monopoles of charge $1$. We collect the main properties of the Prasad--Sommerfield monopole following Atiyah--Hitchin \cite{Atiyah:Hitchin} and Taubes \cite{Jaffe:Taubes}.

First, we fix some notation. For $x \in \R ^{3}$, set $\rho = |x|$ and $\hat{x}=\rho ^{-1}x$. Denote by $\sigma$ the vector $\sigma =(\sigma _{1},\sigma _{2},\sigma _{3}) \in \R ^{3} \otimes \Lie{su}_{2}$, where $\{ \sigma _{1}, \sigma _{2} ,\sigma _{3} \}$ is the standard orthonormal basis of $\Lie{su}_{2}$ defined in terms of Pauli's matrices. In particular, $\sigma _{3} = \diag{(\frac{i}{2},-\frac{i}{2})}$.

The Prasad--Sommerfield (PS) monopole $c_{\Ps}=(A_{\Ps},\Phi _{\Ps})$ is given by the explicit formula, \cf \cite[IV.1, Equation 1.15]{Jaffe:Taubes}:
\begin{alignat}{2}\label{eqn:PS:Monopole}
\Phi _{\Ps}(x)=\left( \frac{1}{\tanh{(\rho )}}-\frac{1}{\rho }\right) \, \hat{x}\cdot \sigma \qquad 
& 
A_{\Ps}(x)=\left( \frac{1}{\sinh{(\rho )}}-\frac{1}{\rho }\right) \, (\hat{x}\times\sigma)\cdot dx
\end{alignat}
Here $\cdot$ and $\times$ are the scalar and vector product in $\R ^{3}$, respectively. To simplify the notation we will often drop the subscript $\Ps$ throughout this section. In \eqref{eqn:PS:Monopole} we fixed the mass $v=1$. A monopole with arbitrary mass $v>0$ is obtained by scaling. The following properties of the PS monopole follow directly from \eqref{eqn:PS:Monopole}, \cf \cite[\S IV.1]{Jaffe:Taubes}.

\begin{lemma}\label{lem:Properties:PS:Monopole}
The pair $(A_{\Ps},\Phi _{\Ps})$ is a solution of the Bogomolny equation (\ref{eqn:Bogomolny}) with finite energy.
\begin{itemize}
\item[(i)] $\Phi$ has exactly one zero, $\Phi (0)=0$.
\item[(ii)] $|\Phi (x)|<1$ and $1-|\Phi (x)|=\frac{1}{\rho} +O(e^{-2\rho})$.
\item[(iii)] By (i), over $\R ^{3}\setminus\{ 0\}$ we can decompose each $\Lie{su}_{2}$--valued form $u$ into diagonal and off-diagonal part $u =u_{D}+u_{T}$, where $u_{D}= |\Phi |^{-2} \langle u,\Phi\rangle\Phi $. Then $|(d_{A}\Phi )_{D}|=O(\rho ^{-2})$ and $|(d_{A}\Phi )_{T}|=O(e ^{-\rho })$.
\end{itemize}
\end{lemma}

In particular, over $\mathbb{R}^{3}\setminus \{ 0\}$ the trivial rank 2 complex vector bundle splits as a sum $H\oplus H^{-1}$ of eigenspaces of $\Phi$. By Lemma \ref{lem:Properties:PS:Monopole}.(i) $H$ is the radial extension of the Hopf line bundle. The adjoint bundle $\left( \R ^{3} \setminus \{ 0 \} \right) \times \Lie{su}_{2}$ splits as a sum $\underline{\mathbb{R}}\oplus H^{2}$. We refer to such a gauge over $\R ^{3} \setminus \{ 0 \}$ as to the \emph{asymptotically abelian gauge} because it yields an asymptotic isomorphism between the PS monopole and a charge $1$ Euclidean Dirac monopole.

The isomorphism $\eta\co \underline{\C}^{2} \ra H\oplus H^{-1}$ can be made explicit, \cf \cite[\S IV.7, 7.1 and 7.2]{Jaffe:Taubes}. A direct computation then yields the following asymptotic expansions.

\begin{lemma}\label{lem:Abelian:Gauge:PS:Monopole}
Let $(A,\Phi )$ be the PS monopole defined in (\ref{eqn:PS:Monopole}).
\begin{itemize}
\item[(i)] There exists an isomorphism $\eta :\underline{\C}^{2} \ra H\oplus H^{-1}$ over $\R ^{3}\setminus\{ 0\}$ such that $\eta (A,\Phi )=(A^{0},\Phi ^{0}) \, \sigma_{3}+(a,\psi )$, with $a$ and $\psi$ a $1$ and $0$--form with values in the $SO(3)$--bundle $\underline{\R} \oplus H^{2}$.
\item[(ii)] $(a,\psi )$ satisfies $d_{ A^{0} \sigma_{3} }^{\ast}a=0=[\Phi ^{0} \sigma_{3},\psi ]$ and $\partial_{\rho}\, \lrcorner\, a=0$. Moreover, as $\rho\rightarrow\infty$:
\[
|a|+|\psi |+|d_{ A^{0}\sigma_{3} }a|+|[\Phi ^{0}\sigma_{3},a]|+|d_{ A^{0}\sigma_{3} }\psi |=O(e^{-\rho })
\]
\end{itemize}
\end{lemma}
Therefore $\eta$ puts $(A,\Phi)$ in ``Coulomb gauge'' with respect to $(A^{0}, \Phi ^{0})\, \sigma _{3}$. Without altering properties (i) and (ii), we have the freedom to change $\eta$ by composing with an element in the stabiliser $U(1)$ of $(A^{0}, \Phi ^{0})\, \sigma _{3}$. By abuse of notation, we won't distinguish between $\eta$ and the induced isomorphism of $SO(3)$--bundles $\left( \R ^{3} \setminus \{ 0 \} \right) \times \Lie{su}_{2} \simeq \underline{ \R } \oplus H^{2}$. As an isomorphism of $SO(3)$--bundles, we have the freedom to compose $\eta$ with an element of $SO(2)$, where $U(1) \ra SO(2)$ is the double cover induced by the adjoint representation $SU(2) \ra SO(3)$.

The deformation theory of charge $1$ monopoles on $\R ^{3}$ can be understood explicitly. Let $D$ be the Dirac operator (\ref{eqn:Dirac:Operator}) twisted by the PS monopole. The $L^{2}$--kernel of $D$ is $4$--dimensional, spanned over $\HH$ by the vector $(d_{A}\Phi ,0)$, \ie
\[
\ker{D}=\langle\, (d_{A}\Phi ,0), \gamma (dx_{h}) \, (d_{A}\Phi ,0 ), h=1,2,3\, \rangle _{\mathbb{R}},
\]
where $\gamma (dx_{h})$ denotes the Clifford multiplication. We can explicitly integrate these infinitesimal deformations. Choose $x_{0}\in\mathbb{R}^{3}$ and let $T_{x_{0}}$ be the translation $x\mapsto x-x_{0}$. Then $T_{x_{0}}^{\ast }c_{\Ps}$ is still a solution to the Bogomolny equation. The corresponding infinitesimal deformation is $-\gamma (x_{0})\, (d_{A}\Phi, 0)=-x_{0} \, \lrcorner \, (F_{A},d_{A}\Phi)$. On the other hand, $(d_{A}\Phi,0)$ is the infinitesimal action of the gauge transformation $\exp{(-\Phi)}$.

For future use, in the next lemma we put $T_{x_{0}}^{\ast }c_{\Ps}$ in ``Coulomb gauge'' with respect to $c_{\Ps}$ and derive some useful estimates.

\begin{lemma}\label{lem:PS:Translations}
There exists $\kappa ,C$ and $\rho _{0}>0$ such that the following holds. For any $x_{0} \in \R ^{3}$ with $|x_{0}| < \kappa$ there exists a solution $(A_{x_{0}},\Phi _{x_{0}})$ to the Bogomolny equation which can be written
\[
(A_{x_{0}},\Phi _{x_{0}}) = (A,\Phi) -x_{0} \, \lrcorner \, (F_{A},d_{A}\Phi) + (a_{x_{0}}, \psi _{x_{0}}), 
\]
where $d_{1}^{\ast}(a_{x_{0}}, \psi _{x_{0}})=0$. Here $d_{1}$ is the linearisation at $c_{\Ps}=(A,\Phi)$ of the action of gauge transformations. Moreover, $|(a_{x_{0}}, \psi _{x_{0}})| \leq C\frac{|x_{0}|^{2}}{\rho ^{3}}$ for all $\rho \geq \rho _{0}$.
\proof
We give a sketch of the proof as it seems that the statement should be well-known.

We will prove later, \cf Proposition \ref{prop:Linearised:Equation:Uj}, that the operator $DD^{\ast}\co W^{2,2}_{w, \delta} \ra L^{2}_{w, \delta -2}$ is an isomorphism for all $\delta \in (-1, 0)$. Here the spaces $W^{m,2}_{w,\delta}$ are defined in Definition \ref{def:Weighted:Spaces:Uj}. Using Lemma \ref{lem:Properties:PS:Monopole}.(v) to see that $d_{A}\Phi \in W^{1,2}_{w,\delta-1}$ one then shows that $DD^{\ast} - 2x_{0}\, \lrcorner \, (F_{A},d_{A}\Phi)\cdot D^{\ast}$ remains an isomorphism if $\kappa$ is sufficiently small. The Implicit Function Theorem implies that there exists a unique solution $u \in W^{2,2}_{w,\delta}$ to the equation
\[
DD^{\ast}u + \big( - x_{0}\, \lrcorner \, (F_{A},d_{A}\Phi) + D^{\ast}u \big) \cdot \big( - x_{0}\, \lrcorner \, (F_{A},d_{A}\Phi) + D^{\ast}u \big) =0,
\]
\ie such that $(A,\Phi) - x_{0}\, \lrcorner \, (F_{A},d_{A}\Phi) + D^{\ast}u $ satisfies the Bogomolny equation and the gauge fixing condition $d_{1}^{\ast}D^{\ast}u=0$. Moreover, $\| D^{\ast}u \| _{W^{1,2}_{w,\delta-1}} \leq C |x_{0}|^{2}$ and the map $x_{0} \mapsto u \in W^{2,2}_{w,\delta}$ is smooth.

It remains to show that $(a_{x_{0}}, \psi _{x_{0}})=D^{\ast}u = O(\rho ^{-3})$. Consider the equation $D\xi + \xi \cdot \xi =0$ for $\xi \in W^{1,2}_{w,\delta -1}$. By Lemma \ref{lem:Abelian:Gauge:PS:Monopole}, in the asymptotically abelian gauge we write
\[
D_{0}\xi + 2 (a,\psi) \cdot  \xi  + \xi \cdot \xi =0,
\]
where $(a,\psi)$ is exponentially decaying and $D_{0}$ is the Dirac operator twisted by the Dirac monopole $(A^{0},\Phi ^{0})$. Decomposing $\xi = \xi_{D}+\xi_T$ into diagonal and transverse part, a Moser iteration argument as in \cite[Lemma 7.10]{Foscolo:Deformation} yields $\xi _{T}=O(\rho ^{-\mu}) = \triangle \xi _{D}$ for all $\mu$.

Apply this result to $\xi = - x_{0}\, \lrcorner \, (F_{A},d_{A}\Phi) + D^{\ast}u$. It follows that we can write
\[
D^{\ast}u=x'_{0} \, \lrcorner \, (F_{A},d_{A}\Phi) + \tau\, (d_{A}\Phi, 0) + O(\rho ^{-3})
\]
for some $x'_{0} \in \R ^{3}$ and $\tau \in \R$. However, since $u \in W^{2,2}_{w,\delta}$ and $\delta \in (-1,0)$, an integration by parts shows that $D^{\ast}u$ is $L^{2}$--orthogonal to $\text{ker}\,{D}$ and therefore $x'_{0}=0=\tau$.
\endproof
\end{lemma}

Of course, there exists a gauge transformation such that $g(A_{x_{0}},\Phi_{x_{0}}) = T^{\ast}_{x_{0}}(A,\Phi)$. Indeed, there must exist $x'_{0} \in \R ^{3}$ and $g$ such that $g(A_{x_{0}},\Phi_{x_{0}}) = T^{\ast}_{x'_{0}}(A,\Phi)$. On the other hand, comparing $|\Phi _{x_{0}}|$ with $|T^{\ast}_{x'_{0}}\Phi|$, one concludes that $x'_{0}=x_{0}$.

\section{A sum of periodic Dirac monopoles}\label{sec:Sum:Dirac:Monopoles}

In this section we construct a singular \emph{reducible} solution $c_{\ext}$ to the Bogomolny equation \eqref{eqn:Bogomolny} by summing periodic Dirac monopoles. Consideration of the asymptotic behaviour of this singular solution suggests the definition of  boundary conditions for periodic monopoles (with singularities) as in Cherkis--Kapustin \cite{Cherkis:Kapustin:1} and \cite{Cherkis:Kapustin:2}. The aim of the rest of the paper will be to construct periodic monopoles by desingularising $c_\ext$ while preserving the boundary conditions.

\medskip
The construction of the pair $c_{\ext} = (A_\ext,\Phi_\ext)$ depends on the choice of the following data:
\begin{itemize}
\item[(i)] A \emph{vacuum background}, \ie constants $v\in\mathbb{R}$ and $b\in\mathbb{R}/\mathbb{Z}$ corresponding to the flat line bundle $L_{v,b}$ on $\RS$ with constant Higgs field $iv$ and flat connection $d+ib\, dt$.
\item[(ii)] The \emph{singularities}: a collection $S$ of $n$ distinct points $p_{i}=(m_{i},a_{i})\in \RS$ for $i=1,\ldots ,n$.
\item[(iii)] The \emph{centres of non-abelian monopoles}: further $k$ points, pair-wise distinct and distinct from the $p_{i}$'s, which we denote by $q_{j}=(z_{j},t_{j})$ for $j=1,\ldots ,k$.
\end{itemize}
Fix an origin in $\mathbb{C} \times \R /2\pi \Z$ so that $z_{1}+\ldots +z_{k}=0$ and $t_{1}+\ldots +t_{k}=0$ modulo $2\pi$. Set $\mu =m_{1}+\ldots +m_{n}$ and $\alpha =a_{1}+\ldots +a_{n}$ modulo $2\pi$. Denote by $d$ the \emph{minimum distance}:
\begin{equation}\label{eqn:Distance}
d=\min{\left\{ |z_{j}-z_{h}|,|z_{j}-m_{i}|\, \text{ for all }\, j,h=1,\ldots k , j\neq h, \text{ and all } i=1,\ldots n\right\}}
\end{equation}
We assume that $d\geq 5$. Throughout the paper constants are allowed to depend on a lower bound for $d$ and on the position of the points $p_{1}, \ldots , p_{n}$.

Define a harmonic function $\Phi_{\ext}$ on $(\RS) \setminus (S \cup \{ q_{1},\ldots ,q_{k}\} )$ by
\begin{equation}\label{eqn:Sum:Dirac:Higgs:Field}
\Phi _{\ext}= v+2\sum _{j=1}^{k}{G _{q_{j}}}-\sum _{i=1}^{n}{G _{p_{i}}}.
\end{equation}
Here $G_{p}$ is the Green's function of $\RS$ with singularity at $p$ given in Lemma \ref {lem:Asymptotics:Periodic:Dirac:Higgs:Field}.

By Lemmas \ref{lem:Asymptotics:Periodic:Dirac:Higgs:Field} and \ref{lem:Asymptotics:Periodic:Dirac:Translations}, for large $|z|$ 
\begin{equation}\label{eqn:Asymptotics:Infinity:Sum:Dirac:Higgs:Field}
\Phi _{\ext}= v +\frac{2k-n}{2\pi}\log{|z|}+\frac{1}{2\pi}\text{Re }\left(\frac{\mu}{z}\right) +O(|z|^{-2}),
\end{equation}
while close to the singularity $p_{i}$
\begin{equation}\label{eqn:Asymptotics:Singularity:Sum:Dirac:Higgs:Field}
\Phi _{\ext}= \text{const} +\frac{1}{2\rho _{i}}+O( \rho _{i} ).
\end{equation}
Here $\rho _{i}=\text{dist}(p_{i},\cdot )$ and the constant term is defined by
\begin{equation}\label{eqn:Mass:Singularity}
\sum _{m \neq i}{ G_{p_{m}}(p_{i}) }-\frac{a_{0}}{2}+v+\frac{1}{\pi}\sum _{j=1}^{k}{\log{|z_{j}-m_{i}|}}+O\left( e^{-d}\right).
\end{equation}
Lemmas \ref{lem:Asymptotics:Periodic:Dirac:Higgs:Field} and \ref{lem:Asymptotics:Periodic:Dirac:Translations} yield similar expansions also for the derivatives of $\Phi _{\ext}$.

Finally, in the ball $B_{\frac{\pi}{2}}(q_{j})$ 
\begin{equation}\label{eqn:Asymptotics:Centre:Higgs:Field:Ext}
\Phi _{\ext}  =  \lambda_{j} -\frac{1}{\rho _{j}} + O\left( \frac{\rho _{j}}{d} + \rho _{j} ^{2} \right)
\end{equation}
where the constant $\lambda _{j}$, using once again Lemma \ref{lem:Asymptotics:Periodic:Dirac:Higgs:Field}, is defined by
\begin{equation}\label{eqn:Scale}
\lambda _{j}=v + a_{0} + \frac{1}{\pi}\sum _{h=1,h\neq j}^{k}{\log{|z_{h}-z_{j}|}}-\frac{1}{2\pi}\sum _{i=1}^{n}{\log{|m_{i}-z_{j}|}}+O\left( e^{-d}\right).
\end{equation}
We will refer to $\lambda _{j}$ as the \emph{local mass} attached to the point $q_{j}$.

\begin{definition}\label{def:Hypothesis:Background:Data}
Given $\lambda _{0},K>1$ and $d_{0} \geq 5$ we say that $v\in\mathbb{R}$ and the points $p_{1},\ldots ,p_{n},q_{1},\ldots ,q_{k}$ are $(\lambda _{0},d_{0},K)$--\emph{admissible} if:
\begin{itemize}
\item[(i)] The minimum distance $d$ (\ref{eqn:Distance}) satisfies $d \geq d_{0}$;
\item[(ii)] $\lambda :=\min _{j}{\lambda _{j}}>\lambda _{0}$;
\item[(iii)] $\overline{\lambda}:=\max _{j}{\lambda _{j}} \leq K \lambda$;
\item[(iv)] $v>0$ if $n=2k$.
\end{itemize}
\end{definition}

The constant $\lambda_0$ has to be thought of as being very large. Therefore it is important to determine which set of data are admissible as $\lambda _{0} \ra \infty$. 
\begin{itemize}
\item[(A)] By \eqref{eqn:Scale}, a first possibility is to fix the points $p_{1}, \ldots ,p_{n}$, $q_{1}, \ldots, q_{k}$ so that (i) in Definition \ref{def:Hypothesis:Background:Data} is satisfied and then choose $v$ sufficiently large. We refer to this as the \emph{high mass} case.
\item[(B)] Consider instead the limit $d \ra +\infty$ and assume that in addition to Definition \ref{def:Hypothesis:Background:Data} we have
\[
\overline{d}=\max{\left\{ |z_{j}-z_{h}|,|z_{j}-m_{i}|\, \text{ for all }\, j,h=1,\ldots k , j\neq h, i=1,\ldots n\right\}} \leq K' d
\]
for some $K'>1$. Then as $d \ra \infty$
\begin{equation}\label{eqn:Limit:Mass:Large:Distance}
\lambda_{j} \sim \frac{1}{\pi}\left( k-1-\frac{n}{2} \right) \log{d}
\end{equation}
for all $j=1, \ldots, k$. Thus if $n<2(k-1)$ we can fix $v$ and $p_{1}, \ldots ,p_{n}$ arbitrarily and then choose the additional $k$ points $q_{1}, \ldots ,q_{k}$ so that $d \geq d_{0}$ for $d_{0}$ sufficiently large. We call this the \emph{large distance} case. It is the most interesting case because, in analogy with Taubes's result \cite{Jaffe:Taubes} for Euclidean monopoles, it is expected that it corresponds to (an open subset of) the end of the moduli space.
\end{itemize}

An important consequence of choosing data with large $\lambda$ is that zeroes of $\Phi _{\ext}$ are localised around the points $q_{1},\dots ,q_{k}$.

\begin{lemma}\label{lem:Localisation:Zeroes:Higgs:Field}
Fix $d_{0}=5$ and suppose that $v > 1$ if $n=2k$. There exists $\lambda _{0}$ such that the following holds. Suppose that the initial data are $(\lambda _{0},5,K)$--admissible. Then for all $j=1, \ldots ,k$ there exists $\frac{1}{\lambda _{j}} < \delta (\lambda _{j}) \leq \frac{2}{\lambda _{j}}$ such that $\Phi _{\ext} \geq 1$ on $(\RS)\setminus \left( S \cup \bigcup _{j=1}^{k}B_{\delta (\lambda_{j} )}(q_{j}) \right)$.
\proof
By Lemmas \ref{lem:Asymptotics:Periodic:Dirac:Higgs:Field} and \ref{lem:Asymptotics:Periodic:Dirac:Translations} there exists a constant $C$ such that in the annulus $\delta \leq \rho _{j}\leq \frac{\pi}{2}$
\[
\Phi _{\ext} \geq \lambda _{j} -\frac{1}{\delta} -C.
\]

Fix some $\lambda _{0}$ such that $\lambda _{0} \geq 1 +\frac{2}{\pi}+C$ and set $\delta = \delta (\lambda_{j} )=\left(\lambda_{j} -1 - C \right) ^{-1}$. By choosing $\lambda _{0}$ larger if necessary, we can assume that $\lambda _{0} \geq 2(1+C)$, so that $\delta (\lambda _{j}) \leq \frac{2}{\lambda _{j}}$.

Now, if $\lambda_{j} >\lambda _{0}$ then $\delta <\frac{\pi}{2}$ and $\Phi _{\ext} \geq 1$ in the annulus $\delta \leq \rho _{j}\leq \frac{\pi}{2}$. On the other hand,  $\Phi _{\ext} \geq 1$ in a neighbourhood of the singularities $p_{1},\ldots ,p_{n}$ and at infinity. Here we need to use the hypothesis $v> 1$ if $n=2k$. The Lemma follows from the minimum principle.
\qed
\end{lemma}

The form $\ast d\Phi _{\ext}$ is the curvature of the line bundle $M \ra (\RS) \setminus \left( S \cup \{ q_{1}, \ldots ,q_{k} \} \right)$
\begin{equation}\label{eqn:Sum:Dirac:Bundle}
M=L_{v,b}\otimes\bigotimes _{j=1}^{k}{L^{2}_{q_{j}}}\otimes \bigotimes _{i=1}^{n}{L^{-1}_{p_{i}}},
\end{equation}
where $L_{q}$ is the line bundle on $(\RS) \setminus \{ p \}$ of Definition \ref{def:Periodic:Dirac:Monopole}.

Consider the reducible $SO(3)$--bundle $\underline{\mathbb{R}}\oplus M$. By a result of Whitney \cite[\S III.7]{Whitney}, isomorphism classes of $SO(3)$--bundles over a CW--complex of dimension at most 3 are completely classified by the second Stiefel--Whitney class $w_{2}$. In our case $w_{2}(\underline{\mathbb{R}}\oplus M)=c_{1}(M)$ modulo $2$, so that $w_{2}$ evaluated on the torus at infinity is $2k-n$ (mod $2$), $1$ on a small sphere enclosing one of the $n$ singularities and it vanishes on spheres enclosing each of the $k$ points.

Denote by $\hat{\sigma}$ the trivialising section of the first factor in $\underline{\R}\oplus M$. Multiplying by $\hat{\sigma}$, $\Phi _{\ext}$ \eqref{eqn:Sum:Dirac:Higgs:Field} defines a Higgs field on $\underline{\R}\oplus M$. Fix the connection
\begin{equation}\label{eqn:Sum:Dirac:Connection}
A_{\ext}=\left( b\, dt+2\sum _{j=1}^{k}{A_{q_{j}}}-\sum _{i=1}^{n}{A_{p_{i}}}\right) \hat{\sigma}
\end{equation}
on $\underline{\R} \oplus M$, where $A_{p}$ is the connection on $L_{p}$ of Definition \ref{def:Periodic:Dirac:Monopole}. Lemmas \ref{lem:Asymptotics:Periodic:Dirac:Connection} and \ref{lem:Asymptotics:Periodic:Dirac:Translations} yield asymptotic expansions of representatives of $A_{\ext}$ close to $p_{i}$, close to $q_{j}$ and at infinity.

\subsection{Boundary conditions}\label{sec:Boundary:Conditions}
By Whitney's result mentioned above, given the collection $S$ of $n$ distinct points $p_{1}, \ldots, p_{n} \in \RS$, there exists a unique $SO(3)$--bundle $V$ up to isomorphism  on $(\RS) \setminus S$ with $w_{2}(V)\cdot [\mathbb{S}^{2}_{p_{i}}]=1$ for all $i=1,\ldots ,n$. Notice that $V$ does not lift to an $SU(2)$--bundle whenever $n>0$ and this explains the choice of $SO(3)$ as structure group. Now, outside of $\{ q_{1}, \ldots, q_{k} \}$ we have an isomorphism $V \simeq \underline{\R} \oplus M$, showing that a desingularisation of the pair $(A_{\ext},\Phi_{\ext})$ is topologically possible. Our aim is to find a desingularisation which solves the Bogomolny equation and satisfies the same asymptotic behaviour of $(A_{\ext},\Phi_{\ext})$ at infinity and close to the points $p_{1}, \ldots, p_{n}$. This leads us to consider the boundary conditions for periodic monopoles (with singularities) introduced by Cherkis and Kapustin in \cite[\S 1.4]{Cherkis:Kapustin:1} and \cite[\S 2]{Cherkis:Kapustin:2}.

\begin{definition}\label{def:Boundary:Conditions}
Given a non-negative integer $k_{\infty} \in \Z _{\geq 0}$, parameters $(v,b) \in \R \times \R / \Z$ and a point $q=(\mu ,\alpha) \in \RS$, let $\mathcal{C}=\mathcal{C}(p_{1},\ldots ,p_{n},k_{\infty},v,b,q)$ be the space of smooth pairs $c=(A,\Phi )$ of a connection $A$ on $V$ and a section $\Phi$ of $V$ satisfying the following boundary conditions.

\begin{enumerate}
\item For each $p_{i} \in S$ there exists a ball $B_{\sigma}(p_{i})$ and a gauge $V|_{B_{\sigma}(p_{i}) \setminus \{ p_{i} \} } \simeq \underline{\R} \oplus H_{p_{i}}$ such that $(A,\Phi)$ can be written
\begin{alignat*}{2}
\Phi = -\frac{1}{2\rho_{i} }\, \hat{\sigma} +\psi\qquad & A = A^{0}\, \hat{\sigma} +a
\end{alignat*}
with $\xi= (a,\psi)=O(\rho_{i}^{-1+\tau})$ and  $|\nabla _{A}\xi| + |[\Phi,\xi]|=O(\rho_{i}^{-2+\tau})$ for some rate $\tau >0$. Here $\rho_{i}$ is the distance from $p_{i}$ and $A^{0}$ is the unique $SO(3)$--invariant connection on $H_{p_{i}}$.

\item There exists $R>0$ and a gauge $V \simeq \underline{\R} \oplus \left( L^{k_{\infty}}_{q} \otimes L_{v,b} \right)$ over $\left( \R^{2} \setminus B_{R} \right) \times \Sph ^{1}$ such that $(A,\Phi)$ can be written
\begin{gather*}
\Phi =\left(v + \frac{k_{\infty}}{2\pi }\log{r} +\frac{1}{2\pi}\text{Re}\left( \frac{\mu}{z}\right)\right) \,\hat{\sigma} +\psi\\
\\
A = \left( b\, dt +k_{\infty}A^{\infty }+\frac{1}{2\pi }(\alpha + k_{\infty}\pi)d\theta +\frac{1}{2\pi}\text{Im}\left( \frac{\mu}{z}\right) dt\right) \, \hat{\sigma} + a
\end{gather*}
with $\xi= (a,\psi)=O(r^{-1-\tau})$ and  $|\nabla _{A}\xi| + |[\Phi,\xi]|=O(r^{-2-\tau})$ for some $\tau >0$. Here $A^{\infty}$ is the connection on $L_{q}$ of Lemma \ref{lem:Asymptotics:Periodic:Dirac:Connection}.
\end{enumerate}
\end{definition}

The reducible pair $(A_{\ext},\Phi_{\ext})$ satisfies these boundary conditions with rate $\tau=1$ and charge at infinity $k_{\infty}=2k-n$. In fact, for topological reasons one always has $k_{\infty} \equiv n$ modulo $2$ and then defines $k$ by the formula above as the (non-abelian) \emph{charge} of the $SO(3)$--monopole $(A,\Phi) \in \mathcal{C}$, \cf \cite[Section 4.2]{Foscolo:Deformation}. We call the parameter $q$ in Definition \ref{def:Boundary:Conditions} the \emph{centre} of the monopole. Observe that it is necessary to fix $q$ in order to have $L^{2}$--integrable infinitesimal deformations.

\section{The family of initial approximate solutions}\label{sec:Approximate:Solutions}

Let $c_{\ext}$ be the reducible solution $(A_{\ext},\Phi _{\ext})$ to the Bogomolny equation obtained in \eqref{eqn:Sum:Dirac:Higgs:Field} and \eqref{eqn:Sum:Dirac:Connection} from $(\lambda _{0},5,K)$--admissible data as in Definition \ref{def:Hypothesis:Background:Data}. Here $\lambda _{0}$ is chosen large enough so that Lemma \ref{lem:Localisation:Zeroes:Higgs:Field} holds. In fact we will reserve the freedom to take $\lambda _{0}$ as large as needed until the end of the construction.

The aim of this section is to construct a family of approximate solutions to the Bogomolny equation \eqref{eqn:Bogomolny}. We desingularise the singular solution $c_{\ext}$ by gluing rescaled PS monopoles in small balls centred at the $k$ points $q_{1}, \ldots, q_{k}$. By varying the centres of the PS monopoles and the gluing maps, we obtain a $(4k-1)$--parameter family of inequivalent smooth pairs $c(x_{0},\tau)$ on $(\RS) \setminus S$.

\begin{definition}\label{def:Gluing:Data}
We define \emph{gluing regions}, adapted cut-off functions and \emph{gluing parameters}.
\begin{itemize}
\item[(i)] Let $N>2$ be a number to be fixed later and set $\delta _{j}=\lambda _{j}^{-\frac{1}{2}}$. Taking $\lambda _{0}$ large enough depending on $N$, we assume that $2N\delta_{j} < \frac{1}{2}$ and $\frac{\delta _{j}}{2N}> \frac{2}{\lambda_{j}}$ for all $j$. Write $(\RS) \setminus S$ as the union of open sets
\begin{alignat*}{2}
U_{j}=B_{N\delta _{j}}(q_{j}) \text{ for } j=1,\ldots ,k, \qquad & U_{\ext}= (\RS) \setminus \left( S \cup \bigcup _{j=1}^{k}{B_{N^{-1}\delta_{j}}(q_{j})}\right).
\end{alignat*}
Let $\Ann_{j}$, $\Ann_{j,\ext}$, $\Ann_{j,\inter}$ be the annuli, respectively,
\[
U_{j}\cap U_{\ext},\qquad B_{2N\delta_{j}}(q_{j})\setminus B_{N\delta_{j}}(q_{j}), \qquad B_{N^{-1}\delta_{j}}(q_{j})\setminus B_{(2N)^{-1}\delta_{j}}(q_{j}).
\]

\item[(ii)] Let $\chi$ be a smooth increasing function of one variable such that $\chi(s) \equiv 1$ if $s \leq 1$ and $\chi (s) \equiv 0$ when $s \geq 2$. For each $j$, let $\rho _{j}$ be the distance function from $q_{j}$ and define cut-off functions $\chi ^{j}_{\inter}=\chi \left( \frac{\delta _{j}\rho_{j}}{2N}\right)$ and $\chi ^{j}_{\ext}=1-\chi \left( N\delta _{j}\rho_{j}\right)$ on $B_{1}(q_{j})$.

\item[(iii)] Fix $\kappa \in (0,1)$ so that Lemma \ref{lem:PS:Translations} holds. Let $\mathcal{P}=\mathcal{P}_{\kappa}$ be the trivial $\T ^{k}$--bundle over the product of $k$ balls $B_{\kappa}(0) \times \dots \times  B_{\kappa}(0) \subset \left( \R ^{3} \right) ^k$. We denote points in $\mathcal{P}$ by $k$--tuples $(x_{0},\tau)$ of points $\left( x^{j}_{0}, \exp{(\tau ^{j}\hat{\sigma})} \right) \in B_{\kappa}(0) \times SO(2)$.
\end{itemize}
\end{definition}

We think of $x^{j}_{0} \in B_{\kappa}(0)$ as parametrising the charge $1$ monopole $(A_{x_{0}^{j}}, \Phi _{x_{0}^{j}})$ on $\mathbb{R}^{3}$ of Lemma \ref{lem:PS:Translations}. Identifying the ball $B_{\pi}(q_{j}) \subset \RS$ with $B_{\lambda _{j}\pi}(0)\subset\mathbb{R}^{3}$ via the homothety 
\begin{equation}\label{eqn:homothety:hj}
h_j\co q_{j}+x \longmapsto \lambda _{j}x,
\end{equation}
the $j$th copy of $B_{\kappa}(0)$ in the definition of $\mathcal{P}$ can also be thought of as $B_{\lambda_{j}^{-1}\kappa}(q_{j}) \subset \RS$.

As for the interpretation of the phase factors $(\tau _{1}, \ldots, \tau _{k})$, let
\begin{equation}\label{eqn:gauge:etaj}
\eta_{j} \co \Ann_{j}\times\mathfrak{su}(2)\rightarrow\underline{\mathbb{R}}\oplus H^{2}\simeq\underline{\mathbb{R}}\oplus M
\end{equation}
be the  isomorphism obtained by composing the gauge transformation $\eta$ of Lemma \ref{lem:Abelian:Gauge:PS:Monopole} with a fixed isomorphism $H^{2} \simeq M$ over $\Ann _{j}$. The choice of $\tau _{j}$ fixes the freedom to compose $\eta _{j}$ with a constant diagonal gauge transformation $\exp{(\tau _{j}\hat{\sigma})}$.

It is clear how to define a family of $SO(3)$--bundles $V(\tau)$ over $(\RS) \setminus S$ with an isomorphism $V(\tau) \simeq \underline{\R} \oplus M$ outside of $\{ q_{1}, \ldots, q_{k} \}$: given a $k$--tuple $\tau = (e^{\tau _{1}\hat{\sigma}}, \ldots, e^{\tau _{k}\hat{\sigma}}) \in SO(2) \times \ldots \times SO(2)$, define $V(\tau)$ identifying $\left( U_{j},U_{j}\times\mathfrak{su}(2)\right)$ and $\left( U_{\ext},\underline{\mathbb{R}}\oplus M\right)$ over $\Ann_{j}$ using $\exp{ (\tau _{j}\hat{\sigma}) } \circ \eta_{j}$. Since $w_{2}\left( V(\tau) \right) \cdot [\Sph ^{2}_{p_{i}}] \equiv 1$ the isomorphism class of $V(\tau)$ does not depend on the choice of $\tau$.

In order to define a smooth pair $c(x_{0},\tau)$ on $V(\tau)$ for all $(x_{0},\tau) \in \mathcal{P}$ we are now going to patch together $c_{\ext}$ and a rescaled Prasad--Sommerfield monopole in a neighbourhood of each of the points $q_{1}, \ldots, q_{k}$. Some care is needed to implement the construction in such a way that the resulting family of approximate solutions to the Bogomolny equation satisfies a number of properties. In particular, the most na\"ive approach to the construction would yield an error term in the Bogomolny equation not sufficiently small to apply the Implicit Function Theorem. Obstructions to match $c_{\ext}$ with scaled PS monopoles at a higher order appear if we also require a fixed behaviour at infinity.

Pull-back a rescaled Prasad--Sommerfield monopole with centre $x^{j}_{0}$ to $U_{j}$ via the homothety $h_{j}$ \eqref{eqn:homothety:hj} obtaining a pair $c_{j}(x_{0})=h_{j}^{\ast}(A_{x_{0}^{j}}, \lambda _{j}\Phi _{x_{0}^{j}})$. By the abelian gauge of Lemma \ref{lem:Abelian:Gauge:PS:Monopole} and Lemma \ref{lem:PS:Translations}, over $U_{j} \setminus \{ q_{j} \}$ we write $e^{\tau_{j}\hat{\sigma}}\eta _{j}\left( c_{j}(x_{0}) \right) = c^{0}_{j}(x_{0}) + (a,\psi)$, where $|(a,\psi)| = O(\lambda _{j}^{-2}\rho_{j} ^{-3})$ and the leading order term $c_{j}^{0}(x_{0})$ is
\begin{equation}\label{eqn:c0:j:x0}
c_{j}^{0}(x_{0}) = c_{j}^{0} -\left( \frac{\langle x \times x^{j}_{0}, dx \rangle}{\lambda _{j}\rho_{j} ^{3}}, \frac{\langle x, x^{j}_{0}\rangle}{\lambda _{j}\rho_{j} ^{3}} \right) \hat{\sigma}.
\end{equation}
Here $c_{j}^{0}$ is the Euclidean Dirac monopole of charge $2$, mass $\lambda _{j}$ and singularity at the origin pulled-back to a neighbourhood of $q_{j}$ in $\RS$. Define a modified pair $c'_{j}(x_{0})$ by
\begin{equation}\label{eqn:c:j:x0:'}
e^{\tau_{j}\hat{\sigma}}\eta _{j}\left( c'_{j}(x_{0}) \right) = c^{0}_{j}(x_{0}) + \chi ^{j}_{\inter}(a,\psi).
\end{equation}

Next, we modify $c_{\ext}$ so that it coincides with $c^{0}_{j}(x_{0})$ over the annulus $\Ann _{j}$. Set
\begin{equation}\label{eqn:c:ext:x0}
c_{\ext}(x_{0})=c_{\ext} -2 \sum_{j=1}^{k}{ \frac{x^{j}_{0}}{ \lambda _{j} } \, \lrcorner \, \left( \ast dG _{q_{j}}, dG_{q_{j}} \right) \otimes \hat{\sigma} }.
\end{equation}
By Lemmas \ref{lem:Asymptotics:Periodic:Dirac:Higgs:Field} and \ref{lem:Asymptotics:Periodic:Dirac:Translations}, over the punctured ball $B_{\frac{\pi}{2}}(q_{j}) \setminus \{ q_{j} \}$ we can write $c_{\ext}(x_{0}) = c_{j}^{0}(x_{0}) + (a,\psi)$ with $|(a,\psi)| = O \left( \frac{\rho _{j}}{d} + \rho _{j}^{2} +\frac{1}{\lambda _{j}}\right)$. Define $c'_{\ext}(x_{0})$ by
\begin{equation}\label{eqn:c:ext:x0:'}
c'_{\ext}(x_{0})  = c^{0}_{j}(x_{0}) + \chi ^{j}_{\ext}(a,\psi).
\end{equation}

The collection $\left( U_{\ext},c'_{\ext}(x_{0})\right)$ and $\left( U_{j},c'_{j}(x_{0}) \right)$ for $j=1, \ldots, k$ defines a pair $c(x_{0},\tau)$ on $V(\tau)$.

\begin{remark}
Notice that the choice $\delta _{j} = \lambda _{j}^{-\frac{1}{2}}$ minimises the size of both $e^{\tau_{j}\hat{\sigma}}\eta _{j}\big( c'_{j}(x_{0}) \big) - c^{j}_{0}(x_{0})$ and $c_{\ext}(x_{0}) - c^{j}_{0}(x_{0})$.
\end{remark}

In fact, as in \cite[Lemma 7.2.46]{Donaldson:Kronheimer}, it is more convenient to fix a base point $\tau_{0}=(\text{id}, \ldots,\text{id})$ and regard the pairs $c(x_{0},\tau)$ as a family of configurations on the fixed $SO(3)$--bundle $V=V(\tau_{0})$. Let $\gamma _{1}, \ldots, \gamma _{k},\gamma _{\ext}$ be a partition of unity subordinate to the cover $U_{1}, \ldots, U_{k},U_{\ext}$. Define a gauge transformation $g_{j}$ on $U_{j}$ by $\eta _{j}\circ g_{j} \circ \eta _{j}^{-1}=\exp{(\tau _{j}\gamma _{\ext}\hat{\sigma})}$. Similarly, let $g_{\ext}$ be the gauge transformation on $U_{\ext}$ with the properties that $g_{\ext} = \exp{(-\tau _{j}\gamma _{j}\hat{\sigma})}$ on $\Ann_{j}$ and $g_{\ext} \equiv 1$ on the complement of $\Ann_{1} \cup \ldots \cup \Ann_{k}$. Then $\eta _{j} \, g_{j}\, \eta _{j}^{-1}\, g_{\ext}^{-1} = e^{\tau_{j}\hat{\sigma}}$ over $\Ann_{j}$ and therefore $\eta _{j}\, g_{j}\left( c'_{j}(x_{0}) \right) = g_{\ext}\left( c'_{\ext}(x_{0}) \right)$. Define a pair on $V$ by
\begin{equation}\label{eqn:c:x0:tau}
\begin{dcases*}
g_{j}\left( c'_{j}(x_{0}) \right) & on $U_{j},$\\
g_{\ext}\left( c'_{\ext}(x_{0}) \right) & on $U_{\ext}.$
\end{dcases*}
\end{equation}
Then $(g_{1}, \ldots, g_{k},g_{\ext})$ defines an isomorphism $g\co V(\tau) \xrightarrow{\sim} V$ such that $g\big( c(x_{0},\tau) \big)$ is precisely \eqref{eqn:c:x0:tau}.

Let $\Gamma$ be the stabiliser of $c_{\ext}$, \ie $\Gamma$ is the group of constant diagonal gauge transformations of $\underline{\R} \oplus M$. There is natural diagonal action of $\Gamma$ on $\tau$ by composition on the left. Since the Prasad--Sommerfield monopole is irreducible $c(x_{0},\tau)$ and $c(x'_{0},\tau ')$ are gauge equivalent if and only if $x'_{0}=x_{0}$ and $\tau '$ belongs to the $\Gamma$--orbit of $\tau$.

\subsection{The centre of mass}\label{sec:Centre:Mass}

By \eqref{eqn:c:ext:x0}, the family $c(x_{0},\tau)$ does not satisfy fixed boundary conditions as $(x_{0},\tau)$ varies in $\mathcal{P}$. Indeed, the centre of the pair $c(x_{0},\tau)$ in the sense of Definition \ref{def:Boundary:Conditions} depends on the \emph{centre of mass}
\begin{equation}\label{eqn:Zeta}
\zeta  = -\sum _{j=1}^{k}{ \frac{ x^{j}_{0} }{\lambda_{j}}}
\end{equation}
of the points $\frac{x_{0}^{1}}{\lambda_{1}}, \ldots, \frac{x_{0}^{k}}{\lambda_{k}}$. Thus the family $\{ c(x_{0},\tau) \, | \, (x_{0},\tau) \in \mathcal{P} \}$ belongs to the fixed configuration space $\mathcal{C}$ of Definition \ref{def:Boundary:Conditions} if and only if $(x_{0},\tau)$ satisfies the ``balancing'' condition $\zeta =0$. Notice that the necessity of this constraint and the action of the stabiliser of $c_{\ext}$ agree with the fact that the dimension of the moduli space of periodic monopoles of charge $k$ is $4(k-1)$ \cite[Theorem 1.5]{Foscolo:Deformation}.

Nonetheless, we will not require $\zeta =0$ at this stage. Since $\RS$ is a \emph{parabolic} manifold, \ie it has no strictly positive Green's function, if $\triangle u = f \in C^{\infty}_{0}(\RS)$, then $u$ grows logarithmically at infinity unless $f$ has mean value zero. As a consequence, when deforming the approximate solution $c(x_{0},\tau)$ into a genuine monopole by the Implicit Function Theorem it is necessary to allow appropriate changes of the asymptotics at infinity by varying the centre $q$ in Definition \ref{def:Boundary:Conditions}. Since our aim is to construct a whole family of solutions to the Bogomolny equation in a fixed moduli space, however, we regard the gluing problem as obstructed. In order to compensate for the obstructions we have to introduce a family of initial approximate solutions depending on parameters. These are precisely the coordinates of the centre of mass $\zeta$. Thus we don't require the ``balancing'' condition $\zeta =0$ at the beginning but will rather fix $\zeta$ only at the end of the construction. This brief discussion motivates the following definition.

\begin{definition}\label{def:Obstruction:Basis}
Define sections $o_{1}, \ldots, o_{4}$ of $(\Lambda ^{1} \oplus \Lambda ^{0}) \otimes V$ over $(\R^2 \times \Sph^1) \setminus S$ by
\begin{alignat*}{2}
o_{h}=-\frac{1}{2\pi k} \sum _{j=1}^{k} { \gamma (dx_{h}) \, \big( \chi ^{j}_{\ext}\, dG_{q_{j}},0 \big) \, \otimes \hat{\sigma}}
\qquad & 
o_{4}=-\frac{1}{2\pi k} \sum _{j=1}^{k} {\big( \chi ^{j}_{\ext}\,  dG_{q_{j}},0 \big) \, \otimes \hat{\sigma} },
\end{alignat*}
Here $h=1,2,3$, $\left( dG_{q_{j}},0 \right)$ is an element of $\Omega(\underline{\R} \oplus M)$ and, by abuse of notation, $dx_{1}=dx$, $dx_{2}=dy$ and $dx_{3}=dt$. Moreover, $\gamma (dx_{h}) \left( a,0 \right)=\partial _{x_{h}} \lrcorner \left( \ast a, a \right)$ is the Clifford multiplication \eqref{eqn:Clifford:Multiplication} of $dx_{h}$ with the $1$--form $a$.
\end{definition}
As we will see later, the span of $o_{1},o_{2},o_{3}$ is the space of obstructions to solve the Bogomolny equation with fixed asymptotics. By Lemmas \ref{lem:Asymptotics:Periodic:Dirac:Higgs:Field} and \ref{lem:Asymptotics:Periodic:Dirac:Translations} and Definition \ref{def:Gluing:Data}.(ii) there exists a constant $C>0$ such that:
\begin{alignat}{2}\label{eqn:Obstruction:Basis}
\begin{dcases*}
|o_{h}| \leq C\rho _{j}^{-2} & in $B_{1}(q_{j}) \setminus B_{N\delta _{j}}(q_{j})$\\
|o_{h}| \leq C & outside of $\bigcup _{j=1}^{k}{ B_{\frac{1}{2}}(q_{j}) }$
\end{dcases*} \qquad &
\begin{dcases*}
|\nabla o_{h}| \leq C\rho _{j}^{-3} & in $B_{1}(q_{j}) \setminus B_{N\delta _{j}}(q_{j})$\\
|\nabla o_{h}| \leq C & outside of $\bigcup _{j=1}^{k}{ B_{\frac{1}{2}}(q_{j}) }$
\end{dcases*}
\end{alignat}

We complete the definition of the pair $c(x_{0},\tau)$ replacing \eqref{eqn:c:x0:tau} by
\begin{equation}\label{eqn:Approximate:Solution}
c(x_{0},\tau)+4\pi \sum_{h=1}^{3} { \zeta_{h}\, o _{h} }
\end{equation}
where $\zeta$ is given by \eqref{eqn:Zeta}. By Lemmas \ref{lem:Asymptotics:Periodic:Dirac:Higgs:Field}, \ref{lem:Asymptotics:Periodic:Dirac:Connection} and \ref{lem:Asymptotics:Periodic:Dirac:Translations} this modification guarantees that the pairs $c(x_{0},\tau)$ satisfy fixed asymptotics for all $(x_{0},\tau) \in \mathcal{P}$, \ie $c(x_{0},\tau)$ lies in the fixed configuration space $\mathcal{C}$ as $(x_{0},\tau)$ varies in $\mathcal{P}$. By abuse of notation, in the rest of the paper we will take \eqref{eqn:Approximate:Solution} as the definition of the pair $c(x_{0},\tau)$.

\subsection{Estimate of the error and the geometry of the approximate solutions.}

Fix a point $(x_{0},\tau) \in \mathcal{P}$ and set $(A,\Phi)=c(x_{0},\tau)$. We want to estimate how far $(A,\Phi)$ is from being a solution to the Bogomolny equation, \ie we want to control $\Psi (x_{0},\tau)=\ast F_{A}-d_{A}\Phi$. For later use, we also estimate the size of the Higgs field $\Phi$ and of the ``curvature'' term $d_{A}\Phi = \ast F_{A}-\Psi(x_{0},\tau)$.

\begin{prop}\label{prop:Pregluing:Map}
There exists $\lambda _{0}$ and $\kappa$ such that the following holds. Suppose that $v,S,q_{1}, \ldots, q_{k}$ are $(\lambda _{0}, 5,K)$--admissible and let $\mathcal{P}=\mathcal{P}_{\kappa}$ be as in Definition \ref{def:Gluing:Data}.(iii). Then the construction of $c(x_{0},\tau)$ in \eqref{eqn:Approximate:Solution} defines a map $c\co \mathcal{P} \ra \mathcal{C}$, the \emph{pre-gluing map}, that factors through $c \co \mathcal{P}/\Gamma \ra \mathcal{C}/\mathcal{G}$, where $\Gamma \simeq SO(2)$ acts on $\mathcal{P}$ diagonally and $\mathcal{G}$ is the space of bounded gauge transformation which preserve the boundary conditions of Definition \ref{def:Boundary:Conditions}.

Furthermore, there exists a uniform constant $C>0$ depending only on $\lambda _{0}$, $\kappa$ and $p_{1}, \ldots, p_{n}$ with the following significance. For every $(x_{0},\tau) \in \mathcal{P}$ define $\zeta$ by \eqref{eqn:Zeta} and set $(A,\Phi) = c(x_{0},\tau)$.
\begin{itemize} 
\item[(i)] The error $\Psi(x_{0},\tau) = \ast F_{A} - d_{A}\Phi$ is supported on $\Ann_{j,\inter}$ and $\Ann_{j,\ext}$. Moreover, define
\begin{equation}\label{eqn:Error:Obstruction}
\Psi _{\zeta} =  4\pi\, d_{2} \left( \sum_{h=1}^{3}{ \zeta_{h}\, o_{h} } \right).
\end{equation}
Then $\Psi _{\zeta}$ is supported on $\Ann_{j,\ext}$ and we have estimates
\begin{alignat*}{3}
|\Psi(x_{0},\tau) - \Psi _{\zeta}| \leq C & \qquad \text{ and } \qquad & \rho _{j}^{2}|\Psi _{\zeta}| \leq \frac{C}{\sqrt{\lambda}}.
\end{alignat*}
\item[(ii)] On every ball $B_{1}(q_{j})$, $j=1, \dots, k$, we have
\begin{equation}\label{eqn:Curvature:Uj}
\left( \lambda _{j}^{-2} + \rho _{j}^{2} \right) |d_{A}\Phi| \leq C.
\end{equation}
\item[(iii)] $|\Phi| \geq \frac{1}{2}$ over $U_{\ext}$.
\end{itemize} 
\proof
The first part of the Proposition is a simple restatement of the construction of the pair $c(x_{0},\tau)$ for $(x_0,\tau) \in \mathcal{P}$. We have to verify the estimates in (i), (ii) and (iii).

From the construction of $c(x_{0},\tau)$ recall that
\[
\begin{dcases*}
c(x_{0},\tau) = c_{j}^{0}(x_{0}) + \chi ^{j}_{\inter}\, O(\lambda _{j}^{-2}\rho _{j}^{-3}) & over $\Ann_{j,\inter}$\\
c(x_{0},\tau) = c_{j}^{0}(x_{0}) + \chi ^{j}_{\ext}\, O(\lambda _{j}^{-1}+\rho _{j}) +  4\pi k \, \sum_{h=1}^{3}{ \zeta_{h}\, o_{h} } & over $\Ann_{j,\ext}$
\end{dcases*}
\]
with $c_{j}^{0}(x_{0})$ defined explicitly in \eqref{eqn:c0:j:x0}. Observe also that if $(A,\Phi )$ and  $(A,\Phi )+(a,\psi )$ both solve the Bogomolny equation and $\chi $ is a smooth function, then
\[
\ast F_{A+\chi a}-d_{A+\chi a}(\Phi +\chi\psi )=\ast (d\chi \wedge a) -(d\chi )\psi +\chi (\chi -1)(\ast [a,a]-[a,\psi ]).
\]
A direct computation using Definition \ref{def:Gluing:Data}.(ii) therefore yields
\begin{equation}\label{eqn:Error}
\begin{dcases*}
| \Psi (x_{0},\tau) | \leq C N^{4} & over $\Ann_{j,\inter}$,\\
|\Psi (x_{0},\tau) - \Psi _{\zeta}| \leq C  & over $\Ann_{j,\ext}$.
\end{dcases*}
\end{equation}
In order to control $|\Psi _{\zeta}|$, we use the fact that $G_{q_{j}}$ is harmonic outside of $q_{j}$ to obtain $|d_{2}o_{h}| \leq C \rho _{j}^{-2}|\nabla \chi ^{j}_{\ext}|$. The estimate now follows from the properties of $\chi ^{j}_{\ext}$ in Definition \ref{def:Gluing:Data}, the definition of $\delta _{j}$ and \eqref{eqn:Zeta}. This concludes the proof of (i).

We prove \eqref{eqn:Curvature:Uj} separately in different regions.
\begin{itemize}
\item[1.] On the ball $\rho _{j} \leq \frac{\delta _{j}}{2N}$, $c(x_{0},\tau)$ is gauge equivalent to a translation of a PS monopole rescaled by $\lambda _{j}$. The quantity $(1+\rho ^{2})|d_{A}\Phi|$ is scale invariant and therefore Lemma \ref{lem:Properties:PS:Monopole} and the fact that $|x^{j}_{0}|<\kappa$ imply the estimate.
\item[2.] On the annulus $\Ann_{j}$, $c(x_{0},\tau)$ coincides with the modified Dirac monopole $c_{j}^{0}(x_{0})$:
\[
\left( \lambda _{j}^{-2} + \rho _{j}^{2} \right) |d_{A}\Phi| \leq C\left( 1+N\lambda _{j}^{-\frac{1}{2}} \right)
\]
follows directly from the definition \eqref{eqn:c0:j:x0} of $c^{j}_{0}(x_{0})$.
\item[3.] We deduce the estimate on the annulus $\Ann_{j,\inter}$ from the previous two and Lemma \ref{lem:PS:Translations}. Write $c(x_{0},\tau) = c^{0}_{j}(x_{0})+ \chi _{j}^{\inter}(a,\psi)$, where $(a,\psi)=O\left( \lambda _{j}^{-2}\rho _{j}^{-3}\right)$. In Lemma \ref{lem:PS:Translations} we did not calculate the decay of the covariant derivative of $(a,\psi)$, but we can argue as follows. Observe that
\[
d_{A+\chi a}\left( \Phi + \chi \psi \right) = d_{A}\Phi + \chi\, d_{A+a}(\Phi + \psi) + d\chi \wedge \psi + \chi (\chi -1)[a,\psi]
\]
where $(A,\Phi)=c^{j}_{0}(x_{0})$ and $\chi = \chi ^{j}_{\inter}$. Since $(A+a,\Phi+\psi)$ is a translation of the PS monopole, we deduce that
\[
\rho _{j}^{2}|d_{A}\Phi| \leq C\big( 1+N^{2}\lambda _{j}^{-1}+ N^{4}\lambda _{j}^{-2} \big).
\]
\item[4.] Finally, on the annulus $B_{1}(q_{j}) \setminus B_{N\delta _{j}}$ write
\[
c(x_{0},\tau)= c_{j}^{0}(x_{0}) + \chi ^{j}_{\ext}(a,\psi) + 4\pi k \sum _{h=1}^{3}{ \zeta _{h}\, o_{h} },
\]
where $(a,\psi)=O(\rho _{j}+\lambda _{j}^{-1})$ and $(\nabla a,\nabla \psi) = O(1)$. A direct computation using \eqref{eqn:Obstruction:Basis} yields $\rho _{j}^{2}|d_{A}\Phi| = O(1+\lambda _{j}^{-\frac{1}{2}})$.
\end{itemize}

Finally, the statement in (iii) is a consequence of the following lemma.
\endproof
\end{prop}

\begin{lemma}\label{lem:Nonvanishing:Higgs:Field}
If $(A,\Phi) = c'_{\ext}(x_{0})+ 4\pi \sum_{h=1}^{3} { \zeta_{h}\, o _{h} }$ then, possibly after taking $\lambda _{0}$ larger and $\kappa$ smaller if necessary, there exists $\frac{1}{\lambda _{j}} < \delta (\lambda _{j}) < \frac{\sqrt{2}}{\lambda _{j}}$ such that
\begin{equation}\label{eqn:Localisation:Zeroes:Higgs:Field}
|\Phi| \geq \frac{1}{2}
\end{equation}
outside of $\bigcup _{j=1}^{k}{ B_{\delta (\lambda _{j})}(q_{j}) }$.
\proof
We prove the lemma controlling the size of $|\Phi|$ through each step of the construction of the pair $c(x_{0},\tau)$. 

First, by the definition \eqref{eqn:c:ext:x0} of $c_{\ext}(x_{0})$, the function $\langle \Phi _{\ext}(x_{0}), \hat{\sigma} \rangle$ is harmonic outside of the points $\{ p_{1}, \ldots, p_{n}, q_{1}, \ldots, q_{k} \}$. One can argue as in Lemma \ref{lem:Localisation:Zeroes:Higgs:Field} to show that there exists $\lambda _{0} > 0$ such that if $\lambda _{j} > \lambda _{0}$ then $\langle \Phi _{\ext}(x_{0}), \hat{\sigma} \rangle \geq 1$ outside of $\bigcup _{j=1}^{k}{ B_{\delta (\lambda _{j})}(q_{j}) }$ for some $\frac{1}{\lambda _{j}}< \delta (\lambda _{j}) \leq \frac{1+\sqrt{1+\kappa}}{2\lambda _{j}}$.

Secondly, picking $\lambda _{0}$ even larger if necessary (depending on $N$), one can make sure that the term of order $O\left( \frac{\rho _{j}}{d}+ \rho _{j}^{2}+ \frac{1}{\lambda _{j}} \right)$ multiplied by the cut-off function $\chi ^{j}_{\ext}$ in \eqref{eqn:c:ext:x0:'} is smaller than $\frac{1}{4}$.

Finally, by \eqref{eqn:Obstruction:Basis}, Definition \ref{def:Hypothesis:Background:Data}.(iii) and the fact that $|\zeta| \leq \frac{C\kappa}{\lambda}$, we can choose $\kappa$ small enough (depending on $K$ of Definition \ref{def:Hypothesis:Background:Data}) so that $\left| 4\pi \sum_{h=1}^{3} { \zeta_{h}\, o_{h} } \right| \leq \frac{1}{4}$.
\endproof
\end{lemma}

\section{The linearised equation}\label{sec:Linear}

Having constructed a family of approximate solutions to the Bogomolny equation, our aim is now to find $\xi=\xi(x_{0},\tau)$ for every $(x_{0},\tau) \in \mathcal{P}$ with the appropriate decay at infinity and at the singularities $p_{i}\in S$ such that $c(x_{0},\tau)+\xi$ is a solution to the Bogomolny equation \eqref{eqn:Bogomolny}. Hence we look for a solution $\xi$ of the equation
\begin{equation}\label{eqn:NonLinear:Equation}
d_{2}\xi + \xi \cdot \xi + \Psi(x_{0},\tau) =0,
\end{equation}
where $d_{2}$ is the linearisation \eqref{eqn:Linearisation:Bogomolny} of the Bogomolny equation at $c(x_{0},\tau)$.

In this section, the technical chore of the paper, we study the linearised equation $d_{2}\xi=f$. This is the crucial step to solve \eqref{eqn:NonLinear:Equation}. The strategy we adopt to understand the invertibility properties of $d_{2}$ is to first solve the equation separately on $U_{j}$ and $U_{\ext}$, in the latter case only modulo obstructions. A global solution of the linearised equation modulo obstructions is then obtained from the local right inverses of $d_{2}$ by a simple iteration. As the main technical tool we will employ a range of weighted Sobolev spaces to carry out the analysis. For technical reasons we will have to distinguish the high mass and large distance limit when studying the equation $d_{2}\xi=f$ over $U_{\ext}$.

Due to the gauge invariance of the Bogomolny equation, $d_{2}$ is not elliptic. We will look for a solution of the form $\xi = d_{2}^{\ast}u$ for a $1$--form $u$ with values in $V$. The Weitzenb\"ock formula
\begin{equation}\label{eqn:Weitzenbock}
d_{2}d_{2}^{\ast}u = \nabla _{A}^{\ast}\nabla _{A}u - \ad^{2}(\Phi)u + \ast[\Psi,u]
\end{equation}
can be deduced from \cite[Lemma 18]{Floer:Monopoles:2}; here $\Psi = \Psi(x_{0},\tau)$. At times it will be convenient to pair a $V$--valued $1$--form $u$ with the zero section of $V$ and consider the $V$--valued form of mixed degree $(u,0)$. Observe that $d_{2}^{\ast}u=D^{\ast}(u,0)$, where $D$ is the Dirac operator \eqref{eqn:Dirac:Operator}, and the equation $d_{2}d_{2}^{\ast}u=f$ is equivalent to
\begin{equation}\label{eqn:d2:D}
DD^{\ast}(u,0)=(f,\ast[\Psi, \ast u]).
\end{equation}
We will make use of the Weitzenb\"ock formulas \cite[Lemma 18]{Floer:Monopoles:2}
\begin{alignat}{3}\label{eqn:Weitzenbock:D}
DD^{\ast}=\nabla _{A}^{\ast}\nabla _{A}-\text{ad}(\Phi)^{2}+\Psi
\qquad & \text{ and } & \qquad
D^{\ast}D=DD^{\ast}+2d_{A}\Phi.
\end{alignat}

\subsection{The linearised equation on $U_{j}$}

There are no obstructions to the invertibility of the operator $d_{2}d_{2}^{\ast}$ over $U_{j}$. The only issue is instead the fact that by \eqref{eqn:Curvature:Uj} the curvature term $d_{A}\Phi$ blows up as $\lambda _{j} \ra \infty$. In view of the Weitzenb\"ock formula \eqref{eqn:Weitzenbock:D} for $D^{\ast}D$, the norm of the inverse of the operator $d_{2}d_{2}^{\ast}\co W^{2,2} \ra L^{2}$ cannot be uniformly bounded. Following a standard approach in gluing problems, we resolve this difficulty introducing appropriate weighted Sobolev spaces.

For each $j=1, \ldots ,k$ consider the weight function $w_{j} = \sqrt{\lambda _{j}^{-2} + \rho _{j}^{2} }$. By abuse of notation, we won't distinguish between the globally defined function $w_{j}$ on $\R ^{3}$ and a fixed smooth increasing function on $(\RS) \setminus S$ with the properties $w_{j} \leq 1$ and
\begin{equation}\label{eqn:Weight:Function:Uj}
w_{j}=\begin{dcases*}
\sqrt{ \lambda _{j}^{-2} + \rho _{j}^{2} } & if $\rho _{j} \leq \frac{1}{2}$,\\
\; 1 & if $\rho _{j} \geq 1$.
\end{dcases*}
\end{equation}

By scaling, we will work on $\R ^{3}$ endowed with the weight function $w=\sqrt{ 1+\rho ^{2} }$. On the trivial $SO(3)$--bundle $\R ^{3} \times \Lie{su}(2)$ we fix a pair $(A,\Phi)$ which coincides with the monopole $(A_{x_{0}},\Phi_{x_{0}})$ of Lemma \ref{lem:PS:Translations} if $\rho \leq (2N)^{-1}\sqrt{\lambda _{j}}$ and with the reducible pair induced by a charge $1$ Euclidean Dirac monopole of unit mass when $\rho \geq N^{-1}\sqrt{\lambda _{j}}$. In other words, we work with the pair obtained from $c'_{j}(x_{0})$ by scaling, but the estimates of Proposition \ref{prop:Linearised:Equation:Uj} below will only depend on the curvature control \eqref{eqn:Curvature:Uj} and the fact that $A$ is a smooth metric connection.

By Kato's inequality, standard results about functions can be extended to $\Lie{su}_{2}$--valued forms and their covariant derivatives. In addition to the Euclidean Sobolev inequality $\| u \| _{L^{6}} \leq C_{Sob}\| \nabla u \| _{L^{2}}$ we will make use of the following Hardy-type inequality, whose proof is obtained by a simple integration by parts \cite[Lemma 2]{Lewis}.
\begin{lemma}\label{lem:Poincare:Inequality:Uj}
For all $\delta \in (-1,0)$ and $u \in C^{\infty}_{0}(\R ^{3}; \Lie{su}_{2})$
\[
\int{ w^{-2\delta -3}\, |u|^{2} } \leq \frac{1}{\delta ^{2}} \int{ w^{-2\delta -1}\, |\nabla _{A}u|^{2} }.
\]
\end{lemma}

\begin{definition}\label{def:Weighted:Spaces:Uj}
For all $\delta \in \R$ and all smooth forms $u \in \Omega (\R ^{3}; \Lie{su}_{2})$ with compact support define:
\begin{alignat*}{2}
\| u \| _{L^{2}_{w,\delta}} = \| w^{-\delta -\frac{3}{2}} u \| _{L^{2}}, \qquad  & \| u \| ^{2}_{W^{1,2}_{w,\delta}} = \| u \| ^{2}_{L^{2}_{w,\delta}} + \| \nabla _{A}u \| ^{2}_{L^{2}_{w,\delta -1}} + \| [\Phi,u] \| ^{2}_{L^{2}_{w,\delta -1}}.
\end{alignat*}
Define spaces $L^{2}_{w,\delta}$ and $W^{1,2}_{w,\delta}$ as the completion of $C^{\infty}_{0}$ with respect to these norms. Finally, we say that $u \in W^{2,2}_{w,\delta}$ if $u \in W^{1,2}_{w,\delta}$ and
\[
\| \nabla _{A}(D^{\ast}u) \| _{L^{2}_{w,\delta -2}} + \| [\Phi, D^{\ast}u] \| _{L^{2}_{w,\delta -2}} < \infty .
\]
\end{definition}

\begin{prop}\label{prop:Linearised:Equation:Uj}
For all $-1 < \delta <0$ there exist $\varepsilon>0$ and $C>0$ such that if $\| w \Psi \| _{L^{3}} < \varepsilon$ then the following holds. For all $f \in L^{2}_{w,\delta -2}$ there exists a unique solution $u \in W^{2,2}_{w,\delta}$ to $d_{2}d_{2}^{\ast}u = f$. Moreover,
\[
\| u \| _{W^{2,2}_{w,\delta}} \leq C \| f \| _{L^{2}_{w,\delta -2}}.
\]
\proof
By an approximation argument, we can assume that $f \in C^{\infty}_{0}$. In view of the Weitzenb\"ock formula \eqref{eqn:Weitzenbock}, the solution $u$ can be found by variational methods. Indeed, by H\"older's inequality
\[
\left| \langle \ast[\Psi,u] , u \rangle _{L^{2}} \right| \leq \| w\Psi \| _{L^{3}} \| w^{-1}u \|_{L^{2}} \| u \| _{L^{6}} \leq 2C_{Sob}\| w\Psi \| _{L^{3}}\| \nabla _{A}u\| _{L^{2}}^{2}.
\]
The last inequality follows from Lemma \ref{lem:Poincare:Inequality:Uj} with $\delta = -\frac{1}{2}$ and the Sobolev inequality. Thus $\| d_{2}^{\ast}u \| _{L^{2}}$ is a norm on $W^{1,2}_{w,\frac{1}{2}}$ provided that $\| w\Psi \| _{L^{3}} < \frac{1}{2C_{Sob}}$. Moreover, since $f\in C^{\infty}_{0}$ the functional $\langle f, u \rangle _{L^{2}}$ is continuous on $W^{1,2}_{w,\frac{1}{2}}$. Hence there exists a unique solution $u$ to $d_{2}d_{2}^{\ast}u = f$, which, by standard elliptic regularity, lies in $W^{1,2}_{w,\frac{1}{2}} \cap C^{\infty}_{loc}$. We have to prove that $u \in W^{2,2}_{w,\delta}$.

In order to justify the integrations by parts below, observe that $|u|=O(\rho ^{-1})$ as $\rho \ra \infty$. This is because $d_{2}d_{2}^{\ast}u=0$ outside of the support of $f$ and $(A,\Phi)$ is reducible outside of a large compact set; thus $u=u_{D}+u_{T}$ with $u_{D}$ harmonic and $u_{T}$ exponentially decaying due to the non-vanishing of $|\Phi|$ at infinity, \cf \cite[Lemma 7.10 and Remark 7.11]{Foscolo:Deformation}. A first integration by parts yields
\begin{align*}
\| u \| _{L^{2}_{w,\delta}} \| f \| _{L^{2}_{w, \delta -2}} &\geq \int{ \langle f, u \rangle\, w^{-2\delta -1} } = \int{ \langle \ast[\Psi,u], u \rangle\, w^{-2\delta -1} }\\
&{} + \int{ (|\nabla _{A}u|^{2} + |[\Phi,u]|^{2}) w^{-2\delta -1}} + (1+2\delta)|\delta | \int{ |u|^{2} w^{2\delta -3} }
\end{align*}
As before, we control the term involving the error $\Psi$ as follows:
\[
\left| \langle \ast[\Psi,u] , u\, w^{-2\delta -1} \rangle _{L^{2}} \right| \leq \| w\Psi \| _{L^{3}} \| u \|_{L^{2}_{w,\delta}} \| w^{-\delta -\frac{1}{2}}u \| _{L^{6}} \leq C  \| w\Psi \| _{L^{3}} \| u \| ^{2}_{W^{1,2}_{w,\delta -1}}
\]
using Lemma \ref{lem:Poincare:Inequality:Uj}, the Sobolev inequality and the fact that $\nabla _{A}(w^{-\delta -\frac{1}{2}}u) \in L^{2}$ if $u \in W^{1,2}_{w,\delta-1}$. Thus
\[
\| u \| _{W^{1,2}_{w,\delta}} \leq C \| f \| _{L^{2}_{w, \delta -2}}
\]
provided that $\| w\Psi \| _{L^{3}}$ is sufficiently small.

Set $\xi = d_{2}^{\ast}u=D^{\ast}(u,0)$. Since $\| \ast[\Psi, \ast u] \| _{L^{2}_{w,\delta -2}} \leq \| w\Psi \| _{L^{3}} \| w^{-\delta -\frac{1}{2}}u \| _{L^{6}}$, the Sobolev inequality and \eqref{eqn:d2:D} imply
 that $\| D\xi \| _{L^{2}_{w,\delta -2}}$ is controlled by $\| f \| _{L^{2}_{w,\delta -2}}$. We will conclude the proof of the Proposition by establishing the a priori estimate
\[
\| \nabla _{A}\xi \| _{L^{2}_{w, \delta -2}} + \| [\Phi, \xi] \| _{L^{2}_{w, \delta -2}} \leq C \left( \| D\xi \| _{L^{2}_{w, \delta -2}} + \| \xi \| _{L^{2}_{w, \delta -1}} \right).
\]
Integrate the Weitzenb\"ock formula \eqref{eqn:Weitzenbock:D} for $D^{\ast}D$ against $w^{-2\delta +1}\xi$:
\begin{align*}
\int{ \left( |\nabla _{A}\xi |^{2}+|[\Phi ,\xi ]|^{2}\right) w^{-2\delta +1} } &\leq c_{1} \int{ |\xi|^{2}w^{-2\delta -1} }+ c_{2} \int{  \langle D\xi ,\xi \rangle\, w^{-2\delta} }\\
&{} + \int{ w^{-2\delta +1} |D\xi|^{2} }+\int{ w^{-2\delta +1} |\Psi | \, |\xi|^{2}}\\
&{} + \int{ w^{-2\delta +1} |d_{A}\Phi| \, |\xi|^{2}},
\end{align*}
using $|\nabla w| \leq 1$. By Proposition \ref{prop:Pregluing:Map}.(ii) $w^{2}|d_{A}\Phi|$ is uniformly bounded. The term involving $\Psi$ is controlled as before using the smallness of $\| w\Psi \| _{L^{3}}$.
\endproof
\end{prop}

\subsection{The linearised equation on $U_{\ext}$: the high mass case}

We move on to study the equation $d_{2}d_{2}^{\ast}u=f$ over the exterior region $U_{\ext}$. The background is now the pair $c'_{\ext}(x_{0}) + 4\pi k \sum_{h=1}^{3} { \zeta_{h}\, o _{h} }$ of \eqref{eqn:c:ext:x0:'} and \eqref{eqn:Approximate:Solution}. Recall that this is a reducible solution to the Bogomolny equation on the complement of $\bigcup _{j=1}^{k}{ B_{2N\delta _{j}}(q_{j}) }$.

If $u$ is a section of the reducible $SO(3)$--bundle $\underline{\R} \oplus M$ of \eqref{eqn:Sum:Dirac:Bundle}, we write $u=u_{D} + u_{T}$. By Lemma \ref{lem:Nonvanishing:Higgs:Field}
\begin{equation}\label{eqn:Control:OffDiagonal:Infinity}
4|[\Phi , u]|^{2} \geq |u_{T}|^{2}
\end{equation}
ouside of small neighbourhoods of $\{ q_{1}, \ldots, q_{k}\}$.
By Fourier analysis with respect to the circle variable $t$ we can further decompose $u_{D} = \Pi _{0}u_{D} + \Pi _{\perp}u_{D}$ into $\Sph ^{1}$--invariant and oscillatory part. On each circle $\{ z \} \times \Sph _{t}^{1}$ one has the Poincar\'e inequality
\begin{equation}\label{eqn:Control:Oscillatory:Infinity}
\int_{\Sph ^{1}}{ |\nabla ( \Pi _{\perp}u_{D} )|^{2} } \geq \int_{\Sph ^{1}}{ |\Pi _{\perp}u_{D}|^{2} }.
\end{equation}

The inequalities \eqref{eqn:Control:OffDiagonal:Infinity} and \eqref{eqn:Control:Oscillatory:Infinity} suggest that, via the Weitzenb\"ock formula \eqref{eqn:Weitzenbock}, we have extremely good control of the off-diagonal and oscillatory piece of $u$ in terms of $d_{2}d_{2}^{\ast}u$. In order to control the $\Sph ^{1}$--invariant diagonal piece $\Pi _{0}u_{D}$ we introduce appropriate weighted spaces. The choice of weight functions is different in the two distinct situations (A) and (B) of Section \ref{sec:Sum:Dirac:Monopoles}, \ie the high mass and large distance case, respectively.
\begin{itemize}
\item[(A)] If we are constructing monopoles in the high mass limit $v \ra \infty$ and $q_{1}, \ldots, q_{k},S$ are contained in a fixed set $B_{R_{0}} \times \Sph ^{1} \subset \RS$, the framework of \cite{Foscolo:Deformation} applies and some care is needed only to check that the constants are uniform as $v \ra \infty$. We will use the weight function $\omega = \sqrt{1+|z|^{2}}$ and let all constants depend on $R_{0}$ without further notice.
\item[(B)] If instead $n \leq 2(k-1)$ and we allow $d \ra \infty$, then the error is concentrated around $k$ points $q_{1}, \ldots, q_{k}$ moving off to infinity and we would like to replace $\sqrt{1+|z|^{2}}$ with a weight function which is uniformly bounded above and below in a neighbourhood of each $q_{j}$ but maintains the same behaviour $O(|z|)$ at large distances.
\end{itemize}

We begin with the high mass case (A) and explain how to extend the results to case (B) in a second step. Thus set $\omega = \sqrt{1+|z|^{2}}$ and consider additional weight functions $\hat{\rho}_{j}, \hat{\rho}_{i}$ defined in a neighbourhood of the point $q_{j}$ and $p_{i} \in S$, respectively. The weight function $\hat{\rho}_{j}$ is a fixed smooth increasing function with $\hat{\rho}_{j} \leq 1$ and
\begin{equation}\label{eqn:Weight:Function:Uext:Singularities}
\hat{\rho}_{j}=\begin{dcases*}
\rho _{j} & if $\rho _{j} \leq \frac{1}{2}$,\\
\; 1 & if $\rho _{j} \geq 1$.
\end{dcases*}
\end{equation}
The function $\hat{\rho}_{i}$ is defined in a similar way, but the transition between $\rho _{i}$ and $1$ takes place on the annulus $B_{2\sigma}(p_{i}) \setminus B_{\sigma}(p_{i})$, where $\sigma >0$ is chosen so that the balls $B_{2\sigma}(p_{i})$ are all disjoint. Constants will be allowed to depend on $\sigma$ without further notice. We denote by $U _{\sigma}$ the complement of the union $\bigcup _{i=1}^{n}{ B_{\sigma}(p_{i}) } \cup \bigcup _{j=1}^{k}{ B_{\frac{1}{2}}(q_{j}) }$. In the definition below, we introduce the relevant Sobolev norms.

\begin{definition}\label{def:Weighted:Spaces:Uext}
Given a triple $(\delta _{1}, \delta _{2}, \delta _{3})$ and a smooth compactly supported form $u \in \Omega (\underline{\R} \oplus M)$ define $\| u \| _{L^{2}_{(\delta _{1}, \delta _{2}, \delta _{3})} }$ as the maximum of the semi-norms:
\begin{alignat*}{3}
\left\| \omega ^{-\delta _{1}-1}u \right\| _{L^{2}(U _{\sigma})}, \qquad & \| \hat{\rho}_{i}^{-\delta _{2}-\frac{3}{2} } u \| _{L^{2}\left( B_{2\sigma}(p_{i}) \right)}, \qquad & \| \hat{\rho}_{j}^{-\delta _{3}-\frac{3}{2} } u \| _{L^{2}\left( B_{1}(q_{j}) \right)}.
\end{alignat*}
Given $\delta >0$, set $\underline{\delta} = (-\delta, \delta , -\delta)$ and for each $m \in \Z$ let $\underline{\delta} + m$ denote the triple $\underline{\delta} + (m,m,m)$. For a smooth compactly supported form $u \in \Omega (\underline{\R} \oplus M)$ we say that
\begin{enumerate}
\item $u \in L^{2}_{\underline{\delta} -2}$ if the corresponding norm is finite;
\item $u \in W^{1,2}_{\underline{\delta} -1}$ if $u \in L^{2}_{\underline{\delta} -1}$ and $\nabla _{A}u, [\Phi, u] \in L^{2}_{\underline{\delta} -2}$;
\item $u \in W^{2,2}_{\underline{\delta}}$ if $D^{\ast}u \in W^{1,2}_{\underline{\delta}-1}$ and $u \in L^{2}_{(\delta, -\delta , -\delta)}$.
\end{enumerate}
The space $W^{m,2}_{\underline{\delta} - 2 +m}$ is defined as the completion of $C^{\infty}_{0}$ with respect to the corresponding norm. By convention $W^{0,2}_{\underline{\delta}-2}=L^{2}_{\underline{\delta}-2}$.
\end{definition}

We stress two aspects of this definition, referring to \cite{Foscolo:Deformation} and the rest of the section for further details. On one side, we distinguished the points $q_{j}$ from the singularities $p_{i}$. More precisely, around each of the singularities $p_{i}$ we imposed a stronger norm. As we will see later, this choice is necessary to control the non-linearities of \eqref{eqn:NonLinear:Equation}. Secondly, observe that if $u \in W^{2,2}_{\underline{\delta}}$ then the transversal and oscillatory part $u_{T}$ and $\Pi _{\perp}u_{D}$ have a stronger decay and lie in $L^{2}_{\underline{\delta}}$. The reason for the odd definition of $W^{2,2}_{\underline{\delta}}$ is to include diagonal sections which have non-zero limits over the punctures and at infinity. This is necessary to ensure the surjectivity of the operator $d_{2}d_{2}^{\ast}$.

The main result of the sub-section is the following proposition.
\begin{prop}\label{prop:Linearised:Equation:Uext}
For all $0< \delta < \frac{1}{2}$ there exists $\varepsilon$ and $C$ with the following significance. Suppose that $\| \hat{\rho}_{j}\Psi |_{B_{1}(q_{j})} \| _{L^{3}} < \varepsilon$ for all $j=1, \ldots ,k$. Then for all $f \in L^{2}_{\underline{\delta}-2}$ such that $\int{ \langle f , \hat{\sigma } \otimes dx_{h} \rangle }=0$ for $h=1,2,3$ there exists a unique solution $\xi \in W^{1,2}_{\underline{\delta}-1}$ to $d_{2}\xi=f$ of the form $\xi =d_{2}^{\ast}u$ with $\int{ \langle u , \hat{\sigma } \otimes dx_{h} \rangle\, \omega ^{-2(\delta +1)} }=0$. Moreover,
\[
\| \xi \| _{W^{1,2}_{\underline{\delta}-1}} \leq C \| f \| _{L^{2}_{\underline{\delta}-2}}.
\]
\end{prop}

The proof is given in two steps. First we prove the existence of a weak solution.

\begin{lemma}\label{lem:Existence:Weak:Solutions:Uext}
For all $0< \delta \leq \frac{1}{2}$ there exists $\varepsilon>0$ and $C$ such that the following holds.

Suppose that $\| \hat{\rho}_{j}\Psi|_{B_{1}(q_{j}) } \| _{L^{3}} < \varepsilon$ for all $j=1, \ldots ,k$. Let $f \in L^{2}_{\underline{\delta}-2}$ be a $(\underline{\R} \oplus M)$--valued $1$--form satisfying $\int{ \langle f , \hat{\sigma } \otimes dx_{h} \rangle }=0$ for $h=1,2,3$. Then there exists a unique  weak solution $u$ to $d_{2}d_{2}^{\ast}u=f$ with
\[
\int{ \langle u , \hat{\sigma } \otimes dx_{h} \rangle\, \omega ^{-2(\delta +1)} }=0 \qquad \text{ and } \qquad \| d_{2}^{\ast}u \| _{L^{2}} \leq C \| f \| _{L^{2}_{\underline{\delta}-2}}.
\]
\proof
The first claim is that $\| d_{2}^{\ast}u \| ^{2}_{L^{2}}$ is uniformly equivalent to $\| \nabla _{A}u \| ^{2}_{L^{2}} + \| [\Phi, u] \| ^{2}_{L^{2}}$ for all $u \in C^{\infty}_{0}$ provided that $\| \hat{\rho}_{j}\Psi|_{B_{1}(q_{j}) } \| _{L^{3}}$ is sufficiently small. Indeed, integrate by parts the Weitzenb\"ock formula \eqref{eqn:Weitzenbock} for $d_{2}d_{2}^{\ast}$ and use the Sobolev and Hardy inequalities to control
\[
\left| \langle \Psi \cdot u , u \rangle _{L^{2}}\right| \leq \| \hat{\rho}_{j}\Psi \| _{L^{3}}\| \hat{\rho}_{j}^{-1}u \| _{L^{2}} \| u \|_{L^{6}} \leq C  \| \hat{\rho}_{j}\Psi \| _{L^{3}} \| \nabla _{A}u \| ^{2}_{L^{2}}.
\]
The second observation is that $\| \nabla _{A}u \| ^{2}_{L^{2}} + \| [\Phi, u] \| ^{2}_{L^{2}}$ is a norm on the space of smooth compactly supported forms with $\int{ \langle u , \hat{\sigma} \rangle \,\omega ^{-2(\delta +1)} } = 0$. More precisely, we are going to show that there exists a uniform constant $C$ such that
\begin{equation}\label{eqn:Existence:Weak:Solutions:Uext}
\| u \| ^{2}_{L^{2}_{( \delta, -\frac{1}{2}, -\frac{1}{2} ) }} \leq C \int{ |\nabla _{A}u|^{2} + |[\Phi, u]|^{2} }.
\end{equation}
Indeed, set $B=B_{1}(q_{j})$ and let $\chi$ be a smooth cut-off function supported in $B$ with $\chi \equiv 1$ in $\frac{1}{2}B$. Lemma \ref{lem:Poincare:Inequality:Uj} applied to $\chi u$ with $\delta = -\frac{1}{2}$ yields
\[
\| \hat{\rho}_{j}^{-1}u \| _{L^{2}(\frac{1}{2}B)} \leq C \left( \| \nabla _{A}u \| _{L^{2}} + \| u \| _{L^{2}(B \setminus \frac{1}{2}B)} \right)
\]
with a uniform constant $C>0$. In the same way we can control the norm $\| \hat{\rho}_{i}^{-1}u\| _{L^{2}}$ in a punctured neighbourhood of $p_{i}$ in terms of $\| \nabla _{A}u \| _{L^{2}}$ and the $L^{2}$--norm of $u$ on an annulus around $p_{i}$. Thus \eqref{eqn:Existence:Weak:Solutions:Uext} will follow once we establish the weighted Poincar\'e inequality
\[
\int{ \omega ^{-2(\delta +1)} |u|^{2} } \leq C \int{ |\nabla _{A}u|^{2} + |[\Phi, u]|^{2} }.
\]
Decompose $u=\Pi_{0}u_{D}+\Pi_{\perp}u_{D}+u_{T}$. Since $\omega \geq 1$, the estimate for $u_{T}$ and $\Pi _{\perp}u_{D}$ follows from \eqref{eqn:Control:OffDiagonal:Infinity} and \eqref{eqn:Control:Oscillatory:Infinity}, respectively. If $u=\Pi_{0}u_{D} \in C^{\infty}_{0}(\R ^{2})$ and $\int{ u \,\omega ^{-2(\delta +1)}} =0$, the stated weighted Poincar\'e inequality is proved in \cite[Corollary 8.4]{Amrouche:Girault:Giroire:1}.

Finally, notice that, for $\delta$ in the range specified, $L^{2}_{\underline{\delta}-2}$ is contained in the dual of $L^{2}_{(\delta, -\frac{1}{2} , -\frac{1}{2} )}$ and therefore the existence of a weak solution $u$ to $d_{2}d_{2}^{\ast}u=f$ follows by variational methods.
\endproof
\end{lemma}

It remains to prove that $\xi \in W^{1,2}_{\underline{\delta}-1}$. This is a consequence of the following a priori estimate.
\begin{lemma}\label{lem:APriori:Estimates:Uext}
For all $0 < \delta <\frac{1}{2}$ there exists $C>0$ such that for all smooth compactly supported $\xi$
\[
\| \xi \| _{W^{1,2}_{\underline{\delta}-1}} \leq C \left( \| D\xi \| _{L^{2}_{\underline{\delta}-2}} + \| \xi \| _{L^{2}} \right).
\]
\proof
The estimate is equation (7.8) in the proof of \cite[Proposition 7.7]{Foscolo:Deformation}. The fact that the constant $C$ is uniform independently of the mass of the monopole is proved in \cite[Lemma 8.11]{Foscolo:Deformation}. For the convenience of the reader we summarise the main aspects of the argument. 

The first step is to prove the weighted elliptic estimate
\begin{equation}\label{eqn:APriori:Estimates:Uext}
\| \xi \| _{W^{1,2}_{\underline{\delta}-1}} \leq C \left( \| D\xi \| _{L^{2}_{\underline{\delta}-2}} + \| \xi \| _{L^{2}_{\underline{\delta}-1}} \right).
\end{equation}
As in the proof of Proposition \ref{prop:Linearised:Equation:Uj}, the main tool is the Weitzenb\"ock formula \eqref{eqn:Weitzenbock:D} for the operator $D^{\ast}D$. The constant $C$ thus depends on appropriate norms of $d_{A}\Phi$ and $\Psi$. Using a partition of unity we can always reduce to prove the estimate under the additional assumption that $\xi$ is supported in a specific given region.
\begin{itemize}
\item[(1)] If $\xi$ is supported in $B_{1}(q_{j})$, as in Proposition \ref{prop:Linearised:Equation:Uj}, integrate by parts the Weitzenb\"ock formula and use Proposition \ref{prop:Pregluing:Map} to show that $|\hat{\rho}_{j}^{2}\Psi|$ and $\hat{\rho}_{j}^{2}|d_{A}\Phi|$ are uniformly bounded.
\item[(2)] Assume now that $\xi \in C^{\infty}_{0}(U_{\sigma})$. Since  $\Psi \equiv 0$ on $U_{\sigma}$, the existence of a uniform constant $C$ follows from the fact that $\omega |d_{A}\Phi|$ is uniformly bounded. In order to check this last statement, recall that $\Phi '_{\ext}(x_{0})$ is a sum of Green's functions and their derivatives. By Lemma \ref{lem:Asymptotics:Periodic:Dirac:Higgs:Field}.(ii), for any $p=(z_{0},t_{0}) \in \RS$
\[
|\nabla G_{p}| + |z-z_{0}|\, |\nabla ^{2}G_{p}| \leq \frac{C}{|z-z_{0}|}
\]
for all $(z,t)$ such that $|z-z_{0}|>2$. We conclude that there exists $C$ depending only on $\sigma$ and $R_{0}$ such that $\omega |d_{A}\Phi| \leq C$.
\item[(3)] Finally, suppose that $\xi$ is supported on $B_{2\sigma}(p_{i})$. In view of the Weitzenb\"ock formula, we have to check that $\hat{\rho }_{i}^{2}|d_{A}\Phi|$ is uniformly bounded. This follows from \eqref{eqn:Asymptotics:Singularity:Sum:Dirac:Higgs:Field} and the fact that the modifications of $\Phi_{\ext}$ in \eqref{eqn:c:ext:x0} and \eqref{eqn:Approximate:Solution} introduce smooth terms controlled by $d^{-1}$ in a neighbourhood of $p_{i}$. 
\end{itemize}

In order to conclude the proof of the Lemma we have to improve \eqref{eqn:APriori:Estimates:Uext} to the required estimate, \ie replace the $L^{2}_{\underline{\delta}-1}$--norm of $\xi$ in the RHS with its $L^{2}$--norm. We distinguish diagonal and off-diagonal part. If $\xi = \xi _{D}$ the statement is deduced from standard theory for the Laplacian in weighted Sobolev spaces, \cf part (1) in the proof of \cite[Proposition 7.7]{Foscolo:Deformation}. Now suppose that $\xi = \xi _{T}$. If $\xi$ is supported on $U_{\sigma}$ we apply word by word the argument of part (2) in the proof of \cite[Proposition 7.7]{Foscolo:Deformation}. Indeed, the argument only relies on \eqref{eqn:Control:OffDiagonal:Infinity} and the fact that $\omega \ra \infty$ as $|z| \ra \infty$. Similarly, the proof of (8.12) in \cite[Lemma 8.11]{Foscolo:Deformation}, case (3), yields the result when $\xi=\xi_{T}$ is supported on $B_{2\sigma}(p_{i})$.
\endproof
\end{lemma}

Proposition \ref{prop:Linearised:Equation:Uext} follows from the two lemmas, since we can control $\| D\xi \| _{L^{2}_{\underline{\delta}-2}}$ in terms of $f$ as in the proof of Proposition \ref{prop:Linearised:Equation:Uj}.

\subsection{The linearised equation on $U_{\ext}$: the large distance case}\label{sec:Linearised:Equation:Uext:Large:Distance}

We come to the task of adapting the analysis to deal with the situation in which the points $q_{1}, \ldots , q_{k}$ move off to infinity. The first step is to define an adapted weight function.

By the assumption $d > d_{0}=5$ of Definition \ref{def:Hypothesis:Background:Data}, the set $B_{2}(z_{j}) \times \mathbb{S}^{1}$ does not contain any of the points $q_{1}, \ldots ,q_{k}, p_{1}, \ldots ,p_{n}$ other than $q_{j}$. By taking $d_{0}$ larger, we can also assume that there exists $R_{0}>0$ such that that the ball $B_{R_{0}}(0) \subset \R ^{2}$ is disjoint from $B_{2}(z_{j})$ for all $j$ and $B_{R_{0}} \times \Sph ^{1}$ contains all the singularities $p_{1}, \ldots, p_{n}$. Set $z_{0}=0$.

In addition to the condition (iii) in Definition \ref{def:Hypothesis:Background:Data}, we will need an extra assumption.

\begin{assumption}\label{assumption}
There exists $K'>1$ such that
\[
\overline{d}=\max{\left\{ |z_{j}-z_{h}|,|z_{j}-m_{i}|\, \text{ for all }\, j,h=1,\ldots k , j\neq h, i=1,\ldots n\right\}} \leq K' d.
\]
\end{assumption}

The assumption clearly implies Definition \ref{def:Hypothesis:Background:Data}.(iii). Moreover, once the centre of mass of $q_{1}, \ldots, q_{k}$ is fixed, this extra requirement is vacuous when $k \leq 2$.

Fix a cover $\{ \Omega_{j} \} _{j=0}^{k}$ of $\RS$ such that $\Omega _{j}$ is an open neighbourhood of the set
\begin{equation}\label{eqn:Voronoi}
\{ (z,t) \in \RS  \text{ such that }  |z-z_{j}| \leq |z-z_{h}| \text{ for all }h=0, \ldots, k \}
\end{equation}
for all $j=0, \ldots, k$. Let $\chi _{0}, \ldots, \chi _{k}$ be a partition of unity subordinate to this cover and set $\omega _{j}(z,t)=\sqrt{1+|z-z _{j}|^{2}}$. We want to define a global weight function $\omega$ such that:
\begin{subequations}\label{eqn:Properties:Weight:Function:Uext}
\begin{align}
&\frac{1}{C_{1}}\, \omega_{j} \leq \omega \leq C_{1}\, \omega _{j} \quad \text{ on }\Omega _{j} \quad \text{ and } \quad \omega \leq C_{1}\omega _{j} \quad\text{ on }\RS \\
&|\nabla \omega| \leq C_{2}, \quad \text{ and }\quad  \left| \omega\,  \triangle \omega \right| \leq C_{3}
\end{align}
\end{subequations}
Given $z_{1}, \ldots ,z_{k} \in \C$, rescale by $d$ around $z_{0}=0$. By Assumption \ref{assumption} $z_{1}, \ldots ,z_{k}$ are mapped to a collection of $k$ points in $\C$ such that the maximum and the minimum of the mutual distances are uniformly bounded above and below. Fix a function $\tilde{r}(z)$ which is a smoothing of $\min_{j= 0,\ldots, k}{\{ |z-d^{-1}z_{j}|\} }$ outside of $z_{0},\ldots, z_{k}$. Since the distance function on $\R ^{2}$ satisfies $|\nabla r| = 1$ and $r\triangle r = -1$ outside of the origin, $\| \nabla \tilde{r} \| _{L^{\infty}}$ and $\| \tilde{r}\triangle \tilde{r}\| _{L^{\infty}}$ are bounded. Now set
\begin{equation}\label{eqn:Weight:Uext:Large:Distance}
\omega(z,t)=\sqrt{1+d^{2}\, \tilde{r}^{2}(d^{-1}z)}.
\end{equation}
Then (\ref{eqn:Properties:Weight:Function:Uext}b) holds with constants depending only on $\| \nabla\tilde{r} \| _{L^{\infty}}$ and $\| \tilde{r}\triangle \tilde{r}\| _{L^{\infty}}$.

Define weighted Sobolev spaces $W^{m,2}_{\underline{\delta} - 2 +m}$, $m=0,1,2$, as in Definition \ref{def:Weighted:Spaces:Uext} using weight functions $\hat{\rho}_{i}, \hat{\rho}_{j}$ and $\omega$. It will also be useful to consider the space $L^{2}_{\underline{\delta}-2,j}$ defined as in Definition \ref{def:Weighted:Spaces:Uext}.(1) using the weight function $\omega _{j}$ instead of $\omega$. By (\ref{eqn:Properties:Weight:Function:Uext}a), $\chi _{j}f \in L^{2}_{\underline{\delta}-2,j}$ for all $f \in L^{2}_{\underline{\delta}-2}$.

We study the equation $d_{2}d_{2}^{\ast}u=f$ in these newly defined spaces. When we restrict to the off-diagonal component, only minor modifications to the proof of Proposition \ref{prop:Linearised:Equation:Uext} are necessary to show that $d_{2}d_{2}^{\ast}\co W^{2,2}_{\underline{\delta}} \ra L^{2}_{\underline{\delta}-2}$ is an isomorphism.

\begin{prop}\label{prop:Linearised:Equation:Uext:Offdiagonal}
For all $0< \delta < \frac{1}{2}$ there exists $\varepsilon$ and $C$ with the following significance. Suppose that $\| \hat{\rho}_{j}\Psi|_{B_{1}(q_{j})} \| _{L^{3}} < \varepsilon$ for all $j=1, \ldots ,k$. Then for all $f=f_{T} \in L^{2}_{\underline{\delta}-2}$ there exists a unique solution $u=u_{T} \in W^{2,2}_{\underline{\delta}}$ to $d_{2}d_{2}^{\ast}u=f$. Moreover,
\[
\| u \| _{W^{2,2}_{\underline{\delta}}} \leq C \| f \| _{L^{2}_{\underline{\delta}-2}}.
\]
\proof
Given that $\omega \geq 1$, $\omega \ra \infty$ as $|z| \ra \infty$ and (\ref{eqn:Properties:Weight:Function:Uext}b) holds, the precise definition of $\omega$ is only used in Lemma \ref{lem:APriori:Estimates:Uext} to show that $\omega |d_{A}\Phi|$ is uniformly bounded on the exterior domain $U_{\sigma}$. Therefore we only have to explain why this quantity remains bounded. Recall that the Higgs field is a sum of Green's function and their derivatives. Then, by Lemma \ref{lem:Asymptotics:Periodic:Dirac:Higgs:Field}.(ii) and (\ref{eqn:Properties:Weight:Function:Uext}a)
\begin{alignat*}{3}
\omega \left( |\nabla G_{q_{j}}| + |\nabla^{2} G_{q_{j}}| \right) \leq C \frac{ \omega _{j} }{|z-z_{j}|} \leq C & \qquad \text{ and } \qquad \omega |dG_{p_{i}}| \leq C \frac{ \omega _{0} }{|z-m_{i}|} \leq C_{i}
\end{alignat*}
if $|z-z_{j}| > 2$ and $|z-m_{i}|>2$, respectively, for a constant $C_{i}$ depending only on $|m_{i}|$.
\qed
\end{prop}

On the diagonal component there is an additional technical difficulty arising from the following finite dimensional family of functions on which the Laplacian is not well-behaved.

\begin{definition}\label{def:vj}
For all $j=0,1,\ldots ,k$ let $\psi_{j}$ be a smooth cut-off function with $\psi _{j} \equiv 0$ if $|z-z_{j}| \leq 1$ and $\psi _{j} \equiv 1$ if $|z-z_{j}| \geq 2$. Define functions $v_{j}=-\frac{1}{4\pi ^{2}}\psi _{j} \log{|z-z_{j}|}$.
\end{definition}
The following two properties of $v_{j}$ are of easy verification.
\begin{itemize}
\item[(i)] There exists a constant $C>0$ such that $\| \nabla v_{j} \| _{L^{\infty}} + \| \nabla ^{2} v_{j} \| _{L^{\infty}} \leq C$;
\item[(ii)] $\int_{\RS}{ \triangle v_{j} }=1$.
\end{itemize}

Given $h \neq j$, set $u=v_{j}-v_{h}$. By (ii) $\triangle u$ has mean value zero and (i) implies that $\| \triangle u \| _{L^{2}_{\underline{\delta}-2} }\leq C$ for a uniform constant $C$. However, restricting to the annulus $2 \leq |z-z_{j}| \leq \frac{1}{2}|z_{j}-z_{h}|$,
\[
\int{ |\nabla u|^{2} } \geq c_{1}\log{|z_{j}-z_{h}|} - \frac{c_{2}}{|z_{j}-z_{h}|^{2}} \xrightarrow{d\ra\infty} \infty
\]
and this fact explains the special role of these functions.

\begin{definition}\label{def:W}
\begin{itemize}
\item[(i)] Let $W$ be the finite dimensional subspace of $1$--forms with values in $\underline{\R} \oplus M$
\[
W= \left\{ \sum_{h=1}^{3}{ \sum_{j=0}^{k}{ \alpha_{h,j}\, d_{2}^{\ast}\left(  v_{j}\, \hat{\sigma} \otimes dx_{h} \right) } } \, \text{ such that } \, \sum _{j=0}^{k}{ \alpha _{h,j} }=0 \text{ for all } h=1,2,3 \right\} .
\]
Define a norm on $W$ by declaring $d_{2}^{\ast}\left(  v_{j}\, \hat{\sigma} \otimes dx_{h} \right)$ an orthonormal system. 
\item[(ii)] Given $f \in L^{2}_{\underline{\delta}-2}$ with $\int{ \langle f, \hat{\sigma} \otimes dx_{h} \rangle }=0$, denote by $\alpha (f)$ the element of $W$ defined by $\alpha _{h,j}=\int{ \langle \chi _{j}f, \hat{\sigma} \otimes dx_{h} \rangle }$.
\end{itemize}
\end{definition}

Notice that the inclusion $L^{2}_{\underline{\delta}-2} \hookrightarrow L^{1}$ is continuous. Indeed, since $\delta >0$, (\ref{eqn:Properties:Weight:Function:Uext}a) implies that $\int{  \omega ^{-2(\delta +1)} } \leq C_{1} \sum _{j=0}^{k}{ \int{  \omega_{j} ^{-2(\delta +1)} } } < +\infty$. In particular, there exists a constant $C>0$ such that
\begin{equation}\label{eqn:Bounded:alpha}
|\alpha (f)| \leq C \| f \| _{L^{2}_{\underline{\delta}-2}}.
\end{equation}

\begin{prop}\label{prop:Linearised:Equation:Uext:Diagonal}
For all $0< \delta < \frac{1}{2}$ there exists a constant $C>0$ with the following significance.

For all $f= f_{D} \in L^{2}_{\underline{\delta}-2}$ such that $\int{ \langle f , \hat{\sigma } \otimes dx_{h} \rangle }=0$ for $h=1,2,3$ there exists $\xi = \xi_{D} \in W^{1,2}_{\underline{\delta}-1}$  such that $d_{2}\xi = f-d_{2}\alpha(f)$. Moreover,
\[
\| \xi \| _{W^{1,2}_{\underline{\delta}-1}} \leq C \| f \| _{L^{2}_{\underline{\delta}-2}}.
\]
\proof
The idea is to write $f=\sum _{j=0}^{k}{\chi _{j}f}$ and use Proposition \ref{prop:Linearised:Equation:Uext} to find a solution to $d_{2}\xi=f$ of the form $\xi = \sum _{j=0}^{k}{d_{2}^{\ast}u_{j}}$.

Set $\alpha _{h,j}=\int { \langle \chi_{j}f, \hat{\sigma} \otimes dx_{h}\rangle }$. The $1$--form $f_{j}=\chi_{j}f-\sum _{h=1}^{3}{\alpha _{h,j}\, d_{2}d_{2}^{\ast} \left( v_{j} \otimes dx_{h}\right) }$ is now orthogonal to constant forms, \ie  $\int{\langle f_{j}, \hat{\sigma} \otimes dx_{h} \rangle } =0$ for all $h=1,2,3$. Moreover, $\| f_{j} \| _{L^{2}_{\underline{\delta}-2,j}} \leq C \| f \| _{L^{2}_{\underline{\delta}-2}}$ by (\ref{eqn:Properties:Weight:Function:Uext}a) and property (i) after Definition \ref{def:vj}. We can therefore apply Proposition \ref{prop:Linearised:Equation:Uext}: there exists $u_{j}$, unique up to the addition of a constant, with the following properties:
\begin{itemize}
\item[(i)] $u_{j}$ is defined on $(\RS) \setminus S$ if $j=0$ and on $(\RS) \setminus \{ q_{j} \}$ otherwise;
\item[(ii)] $d_{2}d_{2}^{\ast} u_{j}=f_{j}$;
\item[(iii)] $\| \omega _{j}^{\delta}d_{2}^{\ast}u _{j} \| _{L^{2}} + \| \omega _{j}^{\delta+1}\nabla (d_{2}^{\ast}u_{j})\| _{L^{2}} \leq C \| f \| _{L^{2}_{\underline{\delta}-2}}$;
\item[(iv)] $ \| u_{0}|_{B_{2\sigma}(p_{i})} \| _{W^{2,2}_{\underline{\delta}}} \leq C \| f \| _{L^{2}_{\underline{\delta}-2}}$ for all $i=1, \ldots, n$;
\item[(v)] $\| u_{j}|_{B_{1}(q_{j})} \| _{W^{2,2}_{\underline{\delta}}} \leq C \| f \| _{L^{2}_{\underline{\delta}-2}}$ for $j \neq 0$.
\end{itemize}

Set $\xi = \sum _{j=0}^{k}{ d_{2}^{\ast}u_{j} }$. We have to show that $\xi \in W^{1,2}_{\underline{\delta}-1}$. Restrict first to the exterior domain $U_{\sigma}$. The fact that $\xi \in W^{1,2}_{\underline{\delta}-1}$ follows immediately from the second statement in (\ref{eqn:Properties:Weight:Function:Uext}a) and (iii) above.

Next, consider the ball $B_{1}(q_{j})$. Since $\omega \geq 1$, (iii) implies that $d_{2}^{\ast}u_{h}|_{B_{1}(q_{j})} \in W^{1,2}$ if $h$ and $j$ are distinct and neither equal to $0$ and similarly $d_{2}^{\ast}u_{h}|_{B_{2\sigma}(p_{i})} \in W^{1,2}$ if $h\neq 0$. Then the Sobolev embedding $W^{1,2} \hookrightarrow L^{6}$ implies that
\[
\| d_{2}^{\ast}u_{h} \| _{W^{1,2}_{\underline{\delta}-1}(B_{1}(q_{j}))} \leq  \| \hat{\rho}_{j}^{\delta -\frac{1}{2}} \| _{L^{3}}\, \| d_{2}^{\ast}u_{h} \| _{L^{6}(B_{1}(q_{j}))} + \| \hat{\rho}_{j}^{\delta +\frac{1}{2}}\| _{L^{\infty}} \, \| \nabla (d_{2}^{\ast}u_{h} ) \| _{L^{2}(B_{1}(q_{j}))} \leq C \| f \| _{L^{2}_{\underline{\delta}-2}}.
\]
The norms of $\hat{\rho}_{j}$ here are bounded because $\delta \in (0,\frac{1}{2})$. Up to changing $\delta$ into $-\delta$, the same argument yields a similar estimate on $B_{2\sigma}(p_{i})$. Together with (iv) and (v) above this concludes the proof.
\qed
\end{prop}

\subsection{Solving the linearised equation modulo obstructions}\label{sec:Linearised:Equation:Modulo:Obstructions}

With these technical details out of the way, we combine Propositions \ref{prop:Linearised:Equation:Uj}, \ref{prop:Linearised:Equation:Uext}, \ref{prop:Linearised:Equation:Uext:Offdiagonal} and \ref{prop:Linearised:Equation:Uext:Diagonal} to solve the equation $d_{2}\xi=f$ modulo obstructions. 

Fix $0< \delta < \frac{1}{2}$. We define weighted Sobolev spaces as in Definition \ref{def:Weighted:Spaces:Uext}, with the difference that over $B_{1}(q_{j})$ we replace $\hat{\rho}_{j}$ of \eqref{eqn:Weight:Function:Uext:Singularities} with the smooth weight function $w_{j}$ of \eqref{eqn:Weight:Function:Uj}. Moreover, the weight function $\omega$ is defined differently in the two situations:
\begin{itemize}
\item[(A)] $S \cup \{ q_{1}, \ldots , q_{k} \} \subset B_{R_{0}} \times \Sph ^{1}$ for some $R_{0}>0$;
\item[(B)] $d \ra \infty$ and  $S, q_{1}, \ldots , q_{k}$ satisfy Assumption \ref{assumption} for some $K'>0$.
\end{itemize}

With these modifications and distinctions understood, the $W^{m,2}_{\underline{\delta}+m-2}$--norm coincides with the $W^{m,2}_{w,\delta}$--norm of Definition \ref{def:Weighted:Spaces:Uj} over $U_{j}$; over $U_{\ext}$, the spaces $W^{m,2}_{\underline{\delta}+m-2}$ are equivalent to the ones used in Proposition \ref{prop:Linearised:Equation:Uext} in case (A) and to those introduced in Propositions \ref{prop:Linearised:Equation:Uext:Offdiagonal} and \ref{prop:Linearised:Equation:Uext:Diagonal} in case (B). In the latter case, consider also the finite dimensional space $W$ introduced in Definition \ref{def:W}.

It is necessary to introduce cut-off functions $\gamma _{j},\gamma _{\ext}, \beta _{j}, \beta _{\ext}$ with some specific properties. Let $\gamma_{j}$ be a smooth function supported in $B_{ 2\delta_{j} }(q_{j})$ and such that $\gamma _{j} \equiv 1$ when $\rho _{j} \leq \frac{\delta_{j}}{2}$. Then $|\nabla \gamma _{j}|\leq \frac{2}{\delta_{j}}$. Define $\gamma _{\ext}$ by $\gamma _{\ext} = 1- \gamma _{j}$ if $\rho \leq 2\delta_{j}$ and $\gamma _{\ext} \equiv 1$ otherwise. The cut-off functions $\beta _{j}$ and $\beta _{\ext}$ are defined in \cite[Lemma 7.2.10]{Donaldson:Kronheimer}: $\beta _{j}$ is a smooth function such that $\beta _{j} \equiv 1$ on $B_{2\delta _{j}}(q_{j})$, $\beta _{j} \equiv 0$ if $\rho _{j} \geq N\delta _{j}$ and
\begin{equation}\label{eqn:CutOff:Gluing}
| \nabla \beta _{j} | \leq \frac{C}{ \log{N} }\frac{1}{\rho _{j}}.
\end{equation}
Similarly, $\beta _{\ext} \equiv 1$ outside of $\bigcup _{j=1}^{k}{ B_{\frac{\delta _{j}}{2}}(q_{j}) }$, $\beta _{\ext} \equiv 0$ on $B_{N^{-1}\delta _{j}}(q_{j})$ and $| \nabla \beta _{\ext} | \leq \frac{C}{ \log{N} }\frac{1}{\rho _{j}}$.

Recall also that in Definition \ref{def:Obstruction:Basis} we distinguished smooth sections $o_{h} \in \Omega(V)$. They are supported on $U_{\ext}$ and under the identification $V|_{U_{\ext}} \simeq \underline{\R} \oplus M$ have only diagonal component. The crucial property of $o_{h}$, $h=1,2,3$, is given by the following lemma.

\begin{lemma}\label{lem:Obstruction:Pairing}
For all $h,l=1,2,3$, $\langle d_{2}o_{h}, \hat{\sigma } \otimes dx_{l} \rangle _{L^{2}} = \delta _{hl}$.
\proof
Since $c(x_{0},\tau)$ is abelian on $U_{\ext}$, $d_{2}o_{h} \oplus d_{1}^{\ast}o_{h}=\slashed{D}o_{h}$, where $\slashed{D}$ is the Dirac operator of $\RS$. Since the Clifford multiplication by $dx_{l}$ commutes with $\slashed{D}$, it is enough to prove that
\begin{alignat*}{2}
\langle \slashed{D}o_{4}, \hat{\sigma } \otimes dx_{h} \rangle _{L^{2}} = 0, & \qquad \langle \slashed{D}o_{4}, \hat{\sigma } \rangle _{L^{2}} = 1.
\end{alignat*}
This follows by direct calculation using the definition $o _{4}=-\frac{1}{2\pi k}\sum_{j=1}^{k}{\left( \chi^{j}_{\ext}\, dG_{q_{j}}, 0 \right) }\hat{\sigma}$. For example,
\[
-\langle \slashed{D}(\chi\, dG,0), \hat{\sigma} \rangle _{L^{2}} = \int_{B} {d(\chi \ast dG)}=\int_{\partial B}{\ast dG}= 2\pi
\]
implies the second identity.
\endproof
\end{lemma}

We refer to $\vspan\{d_{2}o_{h}\, | \, h=1,2,3\} \subset L^{2}_{\underline{\delta}-2}$ as the \emph{obstruction space}. Define $\pi\co L^{2}_{\underline{\delta}-2} \ra L^{2}_{\underline{\delta}-2}$ by
\begin{equation}\label{eqn:Obstruction:Projection}
\pi (f)=f-\sum_{h=1}^{4}{ \langle f , \gamma _{\ext}\, \hat{\sigma }\otimes dx_{h} \rangle _{L^{2}} \,d_{2}o_{h} }.
\end{equation}

\begin{lemma}\label{lem:Continuity:Projection:Obstructions}
There exists a constant $C$ such that
\[
\|\pi(f)\| _{L^{2}_{\underline{\delta}-2}} \leq C \left( N^{-2}\lambda \right) ^{\frac{1-\delta}{2}} \| f \| _{L^{2}_{\underline{\delta}-2}}.
\]
Furthermore, if $f$ is supported on the union of the annuli $\Ann_{j,\inter} \cup \Ann_{j} \cup \Ann_{j,\ext}$ for $j=1, \ldots, k$ then the estimate can be improved to
\[
\| \pi(f) \| _{L^{2}_{\underline{\delta}-2}} \leq C \| f \| _{L^{2}_{\underline{\delta}-2}}.
\]
\proof
Recall that $L^{2}_{\underline{\delta}-2} \hookrightarrow L^{1}$ is continuous. Moreover, if $f$ is supported on $\Ann_{j,\inter} \cup \Ann_{j} \cup \Ann_{j\ext}$
\[
\| f \| _{L^{1}} \leq C \left( N^{-2}\lambda \right) ^{-\frac{1-\delta}{2}} \| f \| _{L^{2}_{\underline{\delta}-2}} 
\]
Therefore, it is enough to estimate $\| d_{2}o_{h}\| _{L^{2}_{\underline{\delta}-2}}$. From Definition \ref{def:Obstruction:Basis} $|d_{2}o_{h}| \leq C \sum _{j=1}^{k}{ \rho ^{-2}_{j}|\nabla \chi ^{j}_{\ext}| }$ and, since $\nabla \chi ^{j}_{\ext}$ is supported in the region where $w_{j} \sim \rho _{j}$ uniformly, we conclude
\begin{equation}\label{eqn:Bound:d2:oh}
\| d_{2}o_{h} \| _{L^{2}_{\underline{\delta}-2}} \leq C \left( N^{-2}\lambda \right) ^{\frac{1-\delta}{2}}. \qedhere
\end{equation}
\end{lemma}

The following theorem yields the solution to the linear problem modulo obstruction.

\begin{thm}\label{thm:Linearised:Equation}
Fix $0< \delta < \frac{1}{2}$.
\begin{itemize}
\item[(A)] Fix $R_{0}>0$ such that $S \cup \{ q_{1}, \ldots, q_{k} \} \subset B_{R_{0}} \times \Sph ^{1}$. There exist $\varepsilon>0$, $N_{0} > 2$ and $C$ with the following significance. Suppose that $\| w_{j}\Psi |_{B_{1}(q_{j})} \| _{L^{3}} < \varepsilon$ for all $j=1, \ldots ,k$ and $N > N_{0}$. Then there exists a map $Q\co \text{im }\pi \subset L^{2}_{\underline{\delta}-2} \ra W^{1,2}_{\underline{\delta}-1}$ such that $\pi \circ d_{2} \circ Q(f)=f$ and
\[
\| Qf \| _{W^{1,2}_{\underline{\delta}-1}} \leq C \| f \| _{L^{2}_{\underline{\delta}-2}}.
\]
\item[(B)] Suppose that $S,q_{1}, \ldots ,q_{k}$ satisfy Assumption \ref{assumption} for some $K'>0$. There exist $\varepsilon>0$, $N_{0} > 2$ and $C>0$ with the following significance. Suppose that $\| w_{j}\Psi|_{B_{1}(q_{j})} \| _{L^{3}} < \varepsilon$ for all $j=1, \ldots ,k$ and $N > N_{0}$. Then there exist a map $Q=(Q_{1},Q_{2})$, where
\[
Q_{1}\co \text{im }\pi \subset L^{2}_{\underline{\delta}-2} \ra W^{1,2}_{\underline{\delta}-1}, \qquad Q_{2} \co \text{im }\pi \subset L^{2}_{\underline{\delta}-2} \ra W,
\]
such that $\pi \circ d_{2} \circ Q_{1}(f)+\pi \circ d_{2}\big( \beta_{\ext}\, Q_{2}(f) \big)=f$. Moreover,
\[
\| Q_{1}(f) \| _{W^{1,2}_{\underline{\delta}-1}} + |Q_{2}(f)| \leq C \| f \| _{L^{2}_{\underline{\delta}-2}}.
\]
\end{itemize}
\proof
We prove the statement in (B). The statement in (A) follows in a similar way using Proposition \ref{prop:Linearised:Equation:Uext} instead of Propositions \ref{prop:Linearised:Equation:Uext:Offdiagonal} and \ref{prop:Linearised:Equation:Uext:Diagonal}.

By abuse of notation, regard $d_{2}$ as the operator $d_{2}\co W^{1,2}_{\underline{\delta}-1} \oplus W \ra L^{2}_{\underline{\delta}-2}$ defined by
\[
d_{2}(\xi,\eta)= d_{2}( \xi + \beta _{\ext}\, \eta ).
\] 
For $f \in L^{2}_{\underline{\delta}-2}$ with $f=\pi(f)$, write $f= \sum _{j=1}^{k}{ \gamma_{j}f }+\gamma_{\ext}f$ and define maps $Q_{1}'\co \text{im }\pi \ra W^{1,2}_{\underline{\delta}-1}$ and $Q_{2}'\co \text{im }\pi \ra W$ as follows: $Q_{2}'(f)=\alpha(\gamma _{\ext}\, f)$, where $\alpha$ is the map of Definition \ref{def:W}.(ii), while
\[
Q_{1}'(f)=\sum_{j=1}^{k}{\beta _{j}\,\xi_{j}} + \beta_{\ext}\, \xi _{\ext}.
\]
Here $\xi_{j}=d_{2}^{\ast}u_{j}$ for the solution $u_{j}$ to $d_{2}d_{2}^{\ast}u_{j}=\gamma _{j}f$ of Proposition \ref{prop:Linearised:Equation:Uj} and $\xi_{\ext}$ is the solution to $d_{2}\xi_{\ext}=\gamma _{\ext}f-\alpha (\gamma _{\ext}f)$ obtained combining Propositions \ref{prop:Linearised:Equation:Uext:Offdiagonal} and \ref{prop:Linearised:Equation:Uext:Diagonal}. In particular, using \eqref{eqn:Bounded:alpha} and \eqref{eqn:CutOff:Gluing}, we deduce the existence of a constant $C>0$ such that
\[
\| Q_{1}'(f) \| _{W^{1,2}_{\underline{\delta}-1}} + |Q_{2} '(f)| \leq C \| f \| _{L^{2}_{\underline{\delta}-2}}.
\]
Set $Q'=(Q_{1}',Q_{2}')$. We have $f-d_{2} Q'(f) = \sum_{j=1}^{k}{ \nabla \beta _{j} \cdot \xi_{j} } + \nabla \beta _{\ext} \cdot \xi_{\ext} + \nabla \beta _{\ext} \cdot \alpha (\gamma _{\ext}f)$. By \eqref{eqn:CutOff:Gluing} and H\"older's inequality
\[
\| f-d_{2} Q'(f)  \| _{L^{2}_{\underline{\delta}-2}} \leq \frac{C}{\log{N} } \| Q_{1}'(f) \| _{W^{1,2}_{\underline{\delta}-1}} + C\frac{\lambda ^{-\frac{1+\delta}{2}}}{\log{N} } \| Q_{2} '(f) \| _{L^{\infty}} \leq \frac{C}{\log{N} }  \| f \| _{L^{2}_{\underline{\delta}-2}}
\] 
because $\| \eta \| _{L^{\infty}} \leq C |\eta|$ for all $\eta \in W$. Since $f-d_{2}  Q'(f) $ is supported in $\bigcup_{j=1}^{k}{\Ann_{j}}$, Lemma \ref {lem:Continuity:Projection:Obstructions} yields $\| f- \pi \circ d_{2}\circ  Q'(f)  \| _{L^{2}_{\underline{\delta}-2}} \leq \frac{C}{\log{N}} \| f \| _{L^{2}_{\underline{\delta}-2}}$ and if $N$ is sufficiently large we can iterate.
\endproof
\end{thm}

\section{Deformation}\label{sec:Deformation}

In this section we complete the construction of a family of solutions to the Bogomolny equation by deforming the approximate solutions $c(x_0,\tau)$ given by the pre-gluing map of Proposition \ref{prop:Pregluing:Map}. Using the projection $\pi$ of \eqref{eqn:Obstruction:Projection}, we split the non-linear equation \eqref{eqn:NonLinear:Equation} into an infinite dimensional and a finite dimensional equation. We solve the infinite dimensional equation first. The following lemma, an immediate consequence of the contraction mapping principle, is a quantitative version of the Implict Function Theorem adapted to case (B) in Theorem \ref{thm:Linearised:Equation}. The statement in case (A) is obtained by setting $W=\{ 0 \}$ (and therefore $\eta =0$).

\begin{lemma}\label{lem:Implicit:Function:Theorem}
Given $c(x_{0},\tau)$ for some $(x_{0},\tau) \in \mathcal{P}$, let $\Psi\co W^{1,2}_{\underline{\delta}-1} \oplus W \rightarrow L^{2}_{\underline{\delta}-2}$ be the smooth map
\[
\Psi (\xi, \eta)=d_{2}(\xi +\beta _{\ext}\eta ) + (\xi +\beta _{\ext}\eta ) \cdot (\xi +\beta _{\ext}\eta ) + \Psi (x_{0},\tau).
\]
Suppose that the following conditions hold.
\begin{itemize}
\item[(i)] There exists a projection $\pi\co  L^{2}_{\underline{\delta}-2} \ra L^{2}_{\underline{\delta}-2}$ such that the map $\pi \circ d_{2}\co W^{1,2}_{\underline{\delta}-1} \oplus W \ra \text{im }\pi$ admits a right inverse $Q=(Q_{1},Q_{2})$ with
\[
\| Q_{1}(f) \| _{ W^{1,2}_{\underline{\delta}-1} } + |Q_{2} (f)| \leq C \| f\| _{L^{2}_{\underline{\delta}-2}}
\]
for all $f \in L^{2}_{\underline{\delta}-2}$ with $\pi(f)=f$.
\item[(ii)] There exists $q>0$ such that
\[
\left\| \pi \Big( (\xi ,\eta ) \cdot (\xi ,\eta )\Big) - \pi\Big( (\xi' ,\eta' ) \cdot (\xi' ,\eta' ) \Big) \right\| _{ L^{2}_{\underline{\delta}-2} } \leq  q \left( \| \xi + \xi' \| _{W^{1,2}_{\underline{\delta}-1} }+| \eta + \eta' |  \right) \left( \| \xi - \xi' \| _{W^{1,2}_{\underline{\delta}-1} }+| \eta - \eta' |  \right) .
\]
Here $(\xi ,\eta ) \cdot (\xi' ,\eta' ) = (\xi +\beta _{\ext}\eta ) \cdot (\xi' +\beta _{\ext}\eta' )$.
\item[(iii)] The error $\Psi (x_{0},\tau)$ satisfies $\| \pi \big( \Psi(x_{0},\tau) \big) \| _{ L^{2}_{\underline{\delta}-2} }\leq \frac{1}{8qC^{2}}$.
\end{itemize}

Then there exists a unique $(\xi ,\eta) \in \text{im }Q \subset W^{1,2}_{\underline{\delta}-1} \oplus W$ such that $\pi \big( \Psi (\xi, \eta) \big)=0$. Moreover,
\[
\| \xi \| _{ W^{1,2}_{\underline{\delta}-1} }+|\eta| \leq 2C\| \pi \big( \Psi (x_{0},\tau) \big) \| _{ L^{2}_{\underline{\delta}-2} }.
\]
\end{lemma}

Theorem \ref{thm:Linearised:Equation} shows that (i) holds provided $\| w_{j}\Psi (x_{0},\tau) \| _{L^{3}} < \varepsilon$ in every ball $B_{1}(q_{j})$, $j=1, \ldots, k$. The next two lemmas imply that (ii) and (iii) are also satisfied if $\lambda$ is sufficiently large.

\begin{lemma}\label{lem:Quadratic:Term}
There exists $C>0$ such that condition (ii) of Lemma \ref{lem:Implicit:Function:Theorem} holds with $q=C\lambda^{\frac{1+\delta}{2}}$.
\proof
Observe that the product $\cdot$ induced by the Clifford multiplication and the Lie bracket is commutative. In particular, $(\xi ,\eta ) \cdot (\xi ,\eta )-(\xi' ,\eta' ) \cdot (\xi' ,\eta' ) = (\xi +\xi', \eta + \eta') \cdot (\xi -\xi', \eta - \eta')$. We will show that $\cdot$ defines a continuous map $\left( W^{1,2}_{\underline{\delta}-1} \oplus W \right) \times \left( W^{1,2}_{\underline{\delta}-1} \oplus W \right)  \ra L^{2}_{\underline{\delta}-2}$ such that
\begin{equation}\label{eqn:Quadratic:Term}
\| (\xi,\eta) \cdot (\xi',\eta') \| _{L^{2}_{\underline{\delta}-2}} \leq C \lambda ^{\delta} \| (\xi, \eta) \| _{W^{1,2}_{\underline{\delta}-1} \oplus W} \| (\xi', \eta') \| _{W^{1,2}_{\underline{\delta}-1} \oplus W}.
\end{equation}
The lemma then follows combining \eqref{eqn:Quadratic:Term} with Lemma \ref{lem:Continuity:Projection:Obstructions}.

The continuity \eqref{eqn:Quadratic:Term} of the product follows from H\"older's inequality and the Sobolev embedding $W^{1,2} \hookrightarrow L^{6}$ along the lines of \cite[Lemmas 5.18 and 6.10]{Foscolo:Deformation}. The crucial observation is that, with respect to the decomposition $V \simeq \underline{\R} \oplus M$ over $U_{\ext}$, there is no $\xi _{D} \cdot \xi'_{D}$ term in the product because $\cdot$ is induced by the Lie bracket on $\Lie{su}(2)$. Moreover, recall that $\beta _{\ext}\eta$ for $\eta \in W$ has only diagonal component. We now provide some details.

Restrict first to the ball $B_{1}(q_{j})$. If $\xi \in W^{1,2}_{\underline{\delta}-1}$ and $\eta \in W$
\[
\int_{B_{1}(q_{j})}{ w_{j}^{2\delta +1}|\xi \cdot \beta_{\ext}\eta|^{2}} \leq C \| w_{j}\eta \|^{2}_{L^{\infty}} \int_{B_{1}(q_{j})}{ w_{j}^{2\delta -1}|\xi |^{2} } \leq C|\eta|^{2} \|\xi \| _{W^{1,2}_{\underline{\delta}-1}}  
\]
because $w_{j} \leq 1$ and $\| \eta \| _{L^{\infty}} \leq C|\eta|$ by Definition \ref{def:W}. On the other hand, in order to estimate the norm $\| \xi \cdot \xi' \| _{L^{2}_{\underline{\delta}-2}}$ for $\xi,\xi' \in W^{1,2}_{\underline{\delta}-1}$, by H\"older's inequality it is enough to observe that
\[
\int_{B_{1}(q_{j})}{ w_{j}^{2\delta +1} |\xi|^{4} } \leq \lambda _{j}^{2\delta}\, \| w_{j}^{\delta -\frac{1}{2}}\xi \| _{L^{2}}\, \| w_{j}^{\delta +\frac{1}{2}}\xi \| ^{3}_{L^{6}} \leq C \lambda ^{2\delta}_{j} \| \xi \| ^{4}_{W^{1,2}_{\underline{\delta}-1}}.
\]
The last inequality follows from the continuity of the Sobolev embedding $W^{1,2}_{\underline{\delta}-1} \hookrightarrow w_{j}^{-\delta -\frac{1}{2}}L^{6}$ on the ball $B_{1}(q_{j})$. The factor of $\lambda _{j}$ in the first inequality is due to the fact that $w_{j} \geq \lambda _{j}^{-1}$.

This establishes \eqref{eqn:Quadratic:Term} on $B_{1}(q_{j})$. The same calculations with $-\delta$ in place of $\delta$ yield the desired estimate also on the ball $B_{2\sigma}(p_{i})$ for all $i=1, \ldots, n$, \cf \cite[Lemmas 5.18]{Foscolo:Deformation}.

On the exterior domain $U _{\sigma }$ write $\xi = \xi _{D} + \xi _{T}$ with respect to the decomposition $V \simeq \underline{\R} \oplus M$. Because of the properties \eqref{eqn:Properties:Weight:Function:Uext} of the weight function $\omega$, if $\xi \in W^{1,2}_{\underline{\delta}-1}$ then $\omega ^{\delta}\xi _{D}, \omega ^{\delta +1}\xi _{T} \in W^{1,2}$. Moreover, $W^{1,2} \hookrightarrow L^{p}$ for all $2 \leq p \leq 6$ by the Sobolev embedding. Therefore
\[
\| \omega ^{\delta +1}(\xi \cdot \xi') \| _{L^{2}} \leq \| \xi \| _{L^{3}} \| \omega ^{\delta +1}\xi '_{T} \| _{L^{6}} + \| \xi' \| _{L^{3}} \| \omega ^{\delta +1}\xi _{T} \| _{L^{6}} \quad \text{ and } \quad \| \omega ^{\delta +1}(\xi \cdot \eta) \| _{L^{2}} \leq \| \eta \| _{L^{\infty}} \| \omega ^{\delta +1}\xi_{T} \| _{L^{2}}
\]
by H\"older's inequality. The estimate \eqref{eqn:Quadratic:Term} follows.
\endproof
\end{lemma}

\begin{remark}
The continuity \eqref{eqn:Quadratic:Term} of the product $\cdot$ justifies the claim that the map $\Psi$ in Lemma \ref{lem:Implicit:Function:Theorem} is smooth.
\end{remark}

\begin{lemma}\label{lem:Smallness:Error}
There exists $C>0$ such that $\| \pi\big( \Psi(x_{0},\tau) \big) \| _{L^{2}_{\underline{\delta}-2}} \leq C\lambda ^{-1-\frac{\delta}{2}}$ and $\| w_{j}\Psi \|_{L^{3}} \leq C \lambda ^{-\frac{1}{2}}$ on each ball $B_{1}(q_{j})$, $j=1, \ldots, k$.
\proof
By Definition \ref{def:Obstruction:Basis} and Lemma \ref{lem:Obstruction:Pairing}, if $\Psi_{\zeta}$ is the component of $\Psi(x_{0},\tau)$ defined by \eqref{eqn:Error:Obstruction} then $\pi (\Psi_{\zeta}) = 0$. By Proposition \ref{prop:Pregluing:Map}.(i)
\begin{equation}\label{eqn:Error:Orthogonal:Obtruction}
\| \Psi(x_{0},\tau) - \Psi_{\zeta} \| _{L^{2}_{\underline{\delta}-2}} \leq C \lambda ^{-1-\frac{\delta}{2}}.
\end{equation}
The second estimate in Lemma \ref{lem:Continuity:Projection:Obstructions} therefore implies the first statement of the lemma.

Similarly, $\| w_{j}\left( \Psi(x_{0},\tau) - \Psi_{\zeta} \right) \| _{L^{3}} = O (\lambda ^{-1})$ on each ball $B_{1}(q_{j})$, while the estimate $\rho _{j}^{2}|\Psi _{\zeta}| \leq \frac{C}{\sqrt{\lambda}}$ in Proposition \ref{prop:Pregluing:Map}.(i) yields $\| w_{j}\Psi _{\zeta} \|_{L^{3}} =O(\lambda ^{-\frac{1}{2}})$.
\endproof
\end{lemma}

\subsection{Existence results}

We have all the ingredients to prove the main result of the paper, an existence theorem for periodic monopoles (with singularities). We begin with a rewriting of the pre-gluing map of Proposition \ref{prop:Pregluing:Map}.

Fix $d_{0}\geq 5$, $K>1$, parameters $v,b$, the set $S$ of singularities and the centres of non-abelian monopoles $q_{1},\dots, q_{k}$. For all $N>2$ we fix $\lambda _{0}(N)$ sufficiently large and assume that $v,S,q_{1}, \ldots , q_{k}$ are $(\lambda _{0},d_{0},K)$--admissible. Moreover, we assume that either:
\begin{itemize}
\item[(A)] There exists $R_{0}>0$ such that $S, \{ q_{1}, \ldots, q_{k} \} \subset B_{R_{0}} \times \Sph ^{1}$; or
\item[(B)] There exist $R_{0}, K'>0$ such that $S \subset B_{R_{0}} \times \Sph ^{1}$, $q_{1}, \ldots, q_{k} \in (\R ^2 \times \Sph^1) \setminus \left( B_{R_{0}} \times \Sph ^{1} \right)$ and Assumption \ref{assumption} is satisfied.
\end{itemize}
Finally, for $\kappa \in (0,1)$ sufficiently small, let $\mathcal{P}=\mathcal{P}_{\kappa}$ be the set of gluing data of Definition \ref{def:Gluing:Data}.

Consider the family $c(x_{0},\tau)$ of Proposition \ref{prop:Pregluing:Map}. If $\kappa$ is sufficiently small and $\lambda _{0}(N)$ sufficiently large, we can assume that \eqref{eqn:Curvature:Uj}, \eqref{eqn:Error} and \eqref{eqn:Localisation:Zeroes:Higgs:Field} are satisfied, uniformly for all $(x_{0},\tau) \in \mathcal{P}$.

\begin{definition}\label{def:Boundary:Conditions:1}
Fix a base point $(0,\tau_{0}) \in \mathcal{P}$ and $\delta >0$ sufficiently small. Set $k_{\infty}=2k-n$ and $q=p_{1}+ \ldots + p_{n}$. Let $\mathcal{C}_{\underline{\delta}} = \mathcal{C}_{\underline{\delta}}(p_{1},\ldots ,p_{n},k_{\infty},v,b,q)$ be the configuration space of pairs $(A,\Phi)$ of the form $(A,\Phi) = c(0,\tau_{0}) + \xi$ with $\xi \in W^{1,2}_{\underline{\delta}-1} \oplus W$. Here and in the rest of the section we set $W=\{ 0 \}$ when condition (A) above holds.
\end{definition}

Observe that once $q_{1}, \ldots, q_{k}$ are fixed the weighted Sobolev norms used to define $\mathcal{C}_{\underline{\delta}}$ are equivalent to those used in \cite{Foscolo:Deformation} to define a smooth structure on the moduli space of periodic monopoles. Then the pregluing map of Proposition \ref{prop:Pregluing:Map} can be considered as a smooth map
\[
c\co \mathcal{P} \ra \mathcal{C}_{\underline{\delta}}.
\]
Smoothness follows from Lemma \ref{lem:PS:Translations} and the explicit construction of $c(x_{0},\tau)$.

The group $\Gamma \simeq SO(2)$ acts on $\mathcal{P}$ by $e^{is} \cdot (x_{0},\tau)=(x_{0},\tau +s)$ and on $\mathcal{C}_{\underline{\delta}}$ as the gauge transformation $\exp{(s\gamma _{\ext}\hat{\sigma})}$; the map $c$ is $\Gamma$--equivariant.

\begin{thm}\label{thm:Existence}
Fix data as above and $N>N_{0}$, where $N_{0}$ is given by Theorem \ref{thm:Linearised:Equation}. Then there exists $\lambda '_{0} \geq \lambda _{0}(N)$ such that if $v,S,q_{1}, \ldots ,q_{k}$ are $(\lambda '_{0},d_{0},K)$--admissible then the following holds.
\begin{itemize}
\item[(i)] There exists a smooth $\Gamma$--equivariant map $c_{1}\co \mathcal{P} \ra \mathcal{C}_{\underline{\delta}}$ such that $\pi \circ \Psi \circ c_{1}=0$. Furthermore, $c_1$ takes the form $c_{1}(x_{0},\tau)=c(x_{0},\tau) + \xi(x_{0},\tau)$ with $\xi(x_{0},\tau) \in W^{1,2}_{\underline{\delta}-1} \oplus W$ and
\[
\| \xi(x_{0},\tau) \| _{W^{1,2}_{\underline{\delta}-1} \oplus W} \leq C\lambda ^{-1-\frac{\delta}{2}}.
\]
\item[(ii)] There exist smooth $\Gamma$--invariant maps $H,h\co \mathcal{P} \ra \R ^{3}$ with
\[
H(x_{0},\tau)=-\sum _{j=0}^{k}{ \frac{x_{0}^{j}}{\lambda _{j}} }
\]
and $|h(x_{0},\tau)| = O(\lambda ^{-\frac{3}{2}})$, such that $c_{1}(x_{0},\tau)$ is a solution to the Bogomolny equation if and only if $H(x_{0},\tau) + h(x_{0},\tau)=0$.
\item[(iii)] Given $(x_{0},\tau) \in \mathcal{P}$ and $\zeta \in \R ^{3}$, let $x_{0}+\zeta$ denote the $k$--tuple $x_{0}+(\zeta, \ldots, \zeta)$. For all $(x_{0},\tau) \in \mathcal{P}_{\frac{\kappa}{2}}$ such that $H(x_{0},\tau)=0$ there exists $\zeta \in \R ^{3}$ such that $|\zeta|=O(\lambda ^{-\frac{1}{2}})$ and
\[
H(x_{0}+\zeta,\tau)+h(x_{0}+\zeta,\tau)=0.
\]
\end{itemize}
\proof
\begin{itemize}
\item[(i)] Consider the map
\[
\pi \circ \Psi\co \mathcal{P} \times \left( W^{1,2}_{\underline{\delta}-1} \oplus W \right) \ra L^{2}_{\underline{\delta}-2},
\]
$(\pi\circ \Psi) (x_{0},\tau, \xi) = \pi \big( d_{2}\xi + \xi\cdot\xi + \Psi(x_{0},\tau) \big)$, which is smooth because of the smoothness of $c\co \mathcal{P} \ra \mathcal{C}_{\underline{\delta}}$, Lemma \ref{lem:Continuity:Projection:Obstructions} and the continuity of the product $\cdot$ established in \eqref{eqn:Quadratic:Term}.

Choosing $\lambda _{0}(N)$ larger if necessary, assume that  $\| w_{j}\Psi(x_{0},\tau) \| _{L^{3}} = O(\lambda ^{-\frac{1}{2}}) < \varepsilon$, where $\varepsilon$ is given by Theorem \ref{thm:Linearised:Equation}. Fix $N>N_{0}$ so that Theorem \ref{thm:Linearised:Equation} holds. Finally, we can choose $\lambda '_{0} \geq \lambda _{0}(N)$ so that
\[
\| \pi\big( \Psi(x_{0},\tau) \big) \| _{L^{2}_{\underline{\delta}-2}} =O(\lambda^{-1-\frac{\delta}{2}}) \leq \frac{1}{2qC^{2}} = O(\lambda^{-\frac{1}{2}-\frac{\delta}{2}})
\]
whenever $\lambda > \lambda '_{0}$. Here $q$ and $C$ are given by Lemma \ref{lem:Quadratic:Term} and Theorem \ref{thm:Linearised:Equation}, respectively, and we used Lemma \ref{lem:Smallness:Error}. Lemma \ref{lem:Implicit:Function:Theorem} therefore implies the existence of the map $c_{1}$. The fact that $c_{1}$ is smooth follows from the fact that the family of right inverses $Q$ of Theorem \ref{thm:Linearised:Equation} depends smoothly on $(x_{0},\tau)$.

\item[(ii)] We are left with the three equations
\[
\langle d_{2}^{\ast}\xi(x_{0},\tau) + \xi(x_{0},\tau) \cdot \xi(x_{0},\tau) +\Psi(x_{0},\tau) , \gamma _{\ext}\,\hat{\sigma }\otimes dx_{h} \rangle_{L^{2}} =0,
\]
for $h=1,2,3$. We have to show that these can be written as $H(x_{0},\tau)+h(x_{0},\tau)=0$ as claimed. Write $\xi(x_{0},\tau)=\xi' + \beta _{\ext}\eta$, with $\xi' \in W^{1,2}_{\underline{\delta}-1}$ and $\eta \in W$. First, we control all the negligible terms:
\begin{itemize}
\item[(a)] $\langle d_{2}\xi', \gamma _{\ext}\,\hat{\sigma }\otimes dx_{h} \rangle_{L^{2}} =O(\lambda ^{-\frac{3}{2}})$. Indeed, integrating by parts
\[
\left| \langle d_{2}\xi', \gamma _{\ext}\,\hat{\sigma }\otimes dx_{h} \rangle_{L^{2}} \right| \leq \| \xi' \| _{ W^{1,2}_{\underline{\delta}-1} }\, \| w_{j}^{-\delta +\frac{1}{2} }\nabla \gamma _{\ext} \| _{L^{2}} = O(\lambda ^{-\frac{3}{2}}).
\]

\item[(b)] $\langle d_{2}(\beta _{\ext}\eta ), \gamma _{\ext}\,\hat{\sigma }\otimes dx_{h} \rangle_{L^{2}} = 0$, because $\gamma _{\ext} \equiv 0$ on the support of $\nabla \beta _{\ext}$, $\beta _{\ext}\equiv 1 \equiv \gamma_{\ext}$ on the support of $d_{2}\eta$ and $\langle d_{2}\eta, \hat{\sigma}\otimes dx_{h} \rangle_{L^{2}} =0$ by  Definition \ref{def:W}.

\item[(c)] By the continuity of the embedding $L^{2}_{\underline{\delta}-2} \hookrightarrow L^{1}$ and \eqref{eqn:Quadratic:Term}
\[
\left| \langle \xi(x_{0},\tau) \cdot \xi(x_{0},\tau) , \gamma _{\ext}\,\hat{\sigma }\otimes dx_{h} \rangle_{L^{2}} \right| \leq C \lambda ^{\delta}\| \xi(x_{0},\tau) \| ^{2}_{W^{1,2}_{\underline{\delta}-1} \oplus W} = O(\lambda^{-2}).
\]

\item[(d)] Finally, as in the proof of Lemma \ref {lem:Continuity:Projection:Obstructions},
\[
\left| \langle \Psi(x_{0},\tau)-\Psi _{\zeta} , \gamma _{\ext} \hat{\sigma}\otimes dx_{h} \rangle_{L^{2}} \right| \leq C \lambda ^{-\frac{1-\delta}{2}} \| \Psi(x_{0},\tau)-\Psi _{\zeta} \| _{L^{2}_{\underline{\delta}-2}} = O(\lambda ^{-\frac{3}{2}})
\]
because $\Psi (x_{0},\tau)-\Psi _{\zeta}$ is supported on $\bigcup_{j=1}^{k}{ (A_{j,\inter}\cup A_{j} \cup A_{j,\ext}) }$.
\end{itemize}
On the other hand, by Lemma \ref {lem:Obstruction:Pairing}
\[
\langle \Psi _{\zeta} , \gamma _{\ext} \hat{\sigma}\otimes dx_{h} \rangle_{L^{2}} = -\sum_{j=1}^{k}{\frac{x^{j}_{0}}{\lambda _{j}}}.
\]

\item[(iii)] The claim follows from Brouwer's Fixed Point Theorem by writing
\[
\zeta = -\frac{h(x_{0}+\zeta,\tau)}{\sum_{j=1}^{k}{\lambda _{j}^{-1}}} = O(\lambda ^{-\frac{1}{2}}).
\]
In the last equality we used Definition \ref{def:Hypothesis:Background:Data}.(iii). \qedhere
\end{itemize}
\end{thm}

In view of \cite[Lemma 7.2]{Foscolo:Deformation}, if the boundary conditions are chosen generically then the monopoles constructed in the theorem are necessary irreducible. Indeed, it is shown in \cite[Lemma 7.2]{Foscolo:Deformation} that there exists no monopole satisfying the boundary conditions of Definition \ref{def:Boundary:Conditions} provided that no subset $\{ p_{i_{1}}, \ldots, p_{i_{k}} \} \subset S$ of cardinality $k$ has centre of mass at the origin. The condition is non-vacuous only when $n \geq k$. In conjunction with \cite[Theorem 1.5]{Foscolo:Deformation}, Theorem \ref{thm:Existence} then shows that, for generic choices of parameters and provided the mass $v$ is sufficiently large when $n \geq 2(k-1)$, the moduli space $\mathcal{M}_{n,k}$ of charge $k$ $SO(3)$ periodic monopoles with $n$ singularities is a non-empty smooth hyperk\"ahler manifold. 

\begin{cor}\label{cor:Existence}
Fix integers $k>0$ and $0 \leq n \leq 2k$, constants $(v,b) \in \R \times \R /\Z$ and distinct points $p_{1}, \ldots, p_{n}, q_{1}, \ldots, q_{k} \in \RS$ with $q_{1}+\ldots+q_{k}=0$. If $n \geq k$ assume that no subset $\{ p_{i_{1}}, \ldots, p_{i_{k}} \}$ has the origin as its centre of mass. Furthermore, assume that conditions (i), (iii) and (iv) of Definition \ref{def:Hypothesis:Background:Data} are satisfied with $d_{0}=5$ and some $K>1$. Let $\mathcal{C}_{\underline{\delta}}=\mathcal{C}_{\underline{\delta}}(p_{1},\ldots ,p_{n},k_{\infty},v,b,q)$ be the configuration space of Definition \ref{def:Boundary:Conditions:1} for some fixed $\delta >0$ sufficiently small.
\begin{itemize}
\item[(i)] Fix $R_{0}>0$ such that $p_{1}, \ldots, p_{n},q_{1}, \ldots, q_{k} \in B_{R_{0}} \times \Sph ^{1}$. Then there exists $v_{0}>0$ such that if $v > v_{0}$ the configuration space $\mathcal{C}_{\underline{\delta}}$ contains irreducible monopoles.
\item[(ii)] Assume that $0 \leq n < 2(k-1)$ and that Assumption \ref{assumption} holds for some $K'>1$. Then there exists $d_{0}>0$ such that if the minimum distance $d$ of \eqref{eqn:Distance} satisfies $d >d_{0}$, then the configuration space $\mathcal{C}_{\underline{\delta}}$ contains irreducible monopoles.
\item[(iii)] If $n=2(k-1)$ and the conditions of (ii) are satisfied, there exist $v_{0}$ and $d_{0}$ such that if $v > v_{0}$ and $d >d_{0}$ then irreducible monopoles exist in $\mathcal{C}_{\underline{\delta}}$.
\end{itemize}
\end{cor}

The statement in (i) establishes the existence of monopoles for any $0 \leq n \leq 2k$ in the high mass case. In case (ii) the mass $v$ is arbitrary and we require that the points $q_{1}, \ldots, q_{k}$ are widely separated. The statement in (iii) follows from \eqref{eqn:Limit:Mass:Large:Distance}: when $n=2(k-1)$, $\lambda _{j}=O(v)$ for all $j$ as $d \ra \infty$. The large distance case considered in (ii) and (iii) is interesting because, in analogy with Taubes's gluing result for Euclidean monopoles, we expect to be able to use Corollary \ref{cor:Existence} to give a description of the asymptotic geometry of the moduli spaces of periodic monopoles (with singularities). We will address this question in a future paper.

\bibliography{Monopoles}
\bibliographystyle{plain}
\end{document}